\documentclass[11pt]{article}
\usepackage{amsmath,amsthm,amssymb}
\usepackage{algorithm,algorithmic}
\usepackage[algo2e,ruled,vlined]{algorithm2e}
\usepackage{enumerate}
\usepackage{hyperref}
\usepackage{tikz}
\usepackage[mathlines]{lineno}
\usepackage{dsfont} 
\usepackage{soul} 
\usepackage{color}

\usepackage{nicefrac} 

\numberwithin{figure}{section}   
\numberwithin{table}{section}   
\numberwithin{equation}{section}  
\numberwithin{algorithm}{section}   

\theoremstyle{plain}
\newtheorem{thm}{Theorem}[section]
\newtheorem{lem}[thm]{Lemma}
\newtheorem{prop}[thm]{Proposition}

\newtheorem{obs}[thm]{Observation}

\newtheorem{quest}[thm]{Question}

\theoremstyle{definition}
\newtheorem{ex}[thm]{Example}
\newtheorem{rem}[thm]{Remark}
\newtheorem{defn}[thm]{Definition}

\definecolor{red}{rgb}{.8,0,0}

\definecolor{bblu}{rgb}{0,0,1}

\definecolor{gray}{rgb}{.4,.4,.4}

\definecolor{gre}{rgb}{0,.7,0}

\newcommand{\mnorm}[1]{\| #1\|_{\rm max}}
\newcommand{\Fnorm}[1]{\| #1\|_{\rm F}}
\newcommand{\nnorm}[1]{\| #1\|_{*}}
\newcommand{\tnorm}[1]{\| #1\|_{2}}
\newcommand{\Snorm}[2]{\| #1\|_{{\rm S},#2}}
\newcommand{\mtx}[1]{\begin{bmatrix} #1 \end{bmatrix}}

\newcommand{\gl}[1]{\operatorname{GL}(#1)}
\newcommand{\SL}[1]{\operatorname{SL}(#1)}
\newcommand{\con}[1]{\operatorname{Con}(#1)}
\newcommand{\nif}[1]{\operatorname{NIF}(#1)}
\newcommand{\Z}{{\mathbb Z}}
\newcommand{\R}{{\mathbb R}}
\newcommand{\C}{{\mathbb C}}

\newcommand{\Cn}{\C^{n}}

\newcommand{\Cnn}{\C^{n\times n}}

\newcommand{\tr}{\operatorname{tr}}

\usepackage[bb=boondox]{mathalfa} 
\newcommand{\ii}{\,\mathbb{i}}
\newcommand{\bone}{\mathbb{1}} 
\newcommand{\bv}{{\bf v}}
\newcommand{\bu}{{\bf u}}
\newcommand{\be}{{\bf e}}
\newcommand{\bw}{{\bf w}}
\newcommand{\bz}{{\bf z}}
\newcommand{\bx}{{\bf x}}
\newcommand{\by}{{\bf y}}

\newcommand{\br}{{\bf r}}
\newcommand{\bzero}{{\bf 0}}

\newcommand{\U}{\mathcal U}
\newcommand{\Q}{\mathcal O}
\newcommand{\lam}{\lambda}

\newcommand{\diag}{\operatorname{diag}}

\newcommand{\spec}{\operatorname{spec}}
\newcommand{\sgn}{\operatorname{sgn}}

\newcommand{\mvec}{\operatorname{vec}}
\newcommand{\spn}{\operatorname{span}}

\newcommand{\ol}{\overline}

\newcommand{\bit}{\begin{itemize}}
\newcommand{\eit}{\end{itemize}}
\newcommand{\ben}{\begin{enumerate}}
\newcommand{\een}{\end{enumerate}}
\newcommand{\beq}{\begin{equation}}
\newcommand{\eeq}{\end{equation}}
\newcommand{\bea}{\begin{eqnarray*}}
\newcommand{\eea}{\end{eqnarray*}}
\newcommand{\bpf}{\begin{proof}}
\newcommand{\epf}{\end{proof}}
\newcommand{\x}{\times}

\newcommand{\lp}{\!\left(}
\newcommand{\rp}{\right)}
\newcommand{\lsb}{\left\{}
\newcommand{\rsb}{\right\}}

\newcommand{\lab}{\left|}
\newcommand{\rab}{\right|}

\DeclareMathOperator*{\rank}{rank}

\begin{document}

\title{Apportionable matrices and gracefully labelled graphs}
\author{Antwan Clark\thanks{Department of Applied Mathematics and Statistics, Johns Hopkins University, Baltimore, MD 21218, USA (aclark1@jhu.edu)}\and Bryan A.~Curtis\thanks{Department of Mathematics, Iowa State University,
Ames, IA 50011, USA (bcurtis1@iastate.edu)}\and Edinah K.~Gnang\thanks{Department of Applied Mathematics and Statistics, Johns Hopkins University, Baltimore, MD 21218, USA (egnang1@jhu.edu)}\and Leslie Hogben\thanks{Department of Mathematics, Iowa State University,
Ames, IA 50011, USA and  American Institute of Mathematics, Caltech, 1200 E California Blvd, Pasadena, CA 91125, USA
(hogben@aimath.org)}} 

\maketitle\vspace{-15pt}

\begin{abstract} 
To apportion a complex matrix means to apply a similarity so that all entries of the resulting matrix have the same magnitude.  We initiate the study of apportionment, both by unitary matrix similarity and general matrix similarity. There are connections between apportionment and classical  graph decomposition problems, including graceful labelings of graphs, Hadamard matrices, and equiangluar lines, and potential applications to instantaneous uniform mixing in quantum walks.  The connection between apportionment and graceful labelings allows the construction of apportionable matrices from trees. A generalization of the well-known Eigenvalue Interlacing Inequalities using graceful labelings is also presented. It is shown that every rank one matrix can be apportioned by a unitary similarity, but there are $2\x 2$ matrices that cannot be apportioned.  A necessary condition for a matrix to be apportioned by  unitary matrix is established. This condition is used to construct a set of matrices with nonzero Lebesgue measure that are not apportionable by a unitary matrix. 
\end{abstract}

\noindent\textbf{Keywords.} uniform matrix, apportionment, unitary, graceful labeling

\noindent\textbf{AMS subject classifications.} 05C50, 05C78, 15A04, 15A18

\section{Introduction}\label{s:intro}

There has been extensive study of diagonalization of matrices,  and finding the Jordan canonical form for  matrices that  are not diagonalizable. Diagonalization can be viewed as using a similarity to concentrate the magnitude of all the entries within a small subset of entries. Here we study what can be viewed as reversing this process, spreading out the magnitudes as uniformly as possible. A square complex matrix is \emph{uniform} if all entries have the same absolute value. Hadamard matrices and discrete Fourier transforms are important examples of uniform matrices.  Uniform matrices play the role of the target (like diagonal matrices) in this process. A square complex matrix is \emph{apportionable} if it is similar (by a specified type of  matrix) to a uniform matrix (formal definitions of various types of apportionability and other terms are given below).  There are interesting connections between apportionability and classical combinatorial problems, such as graceful labeling of graphs (see Section \ref{s:graph-decomp}), construction of equiangular lines  (see Section \ref{s:nec-translate}), and construction of Hadamard matrices  (see discussion of these matrices later in this section).  There are also connections with the relatively new study of instantaneous uniform mixing in continuous-time quantum walks. 
Specifically, the continuous-time quantum walk on a simple  graph $G$ has the transition operator $e^{-itA_G}$, where $A_G$ is the adjacency matrix of $G$, and has instantaneous uniform mixing at time $t_0$ if and only if $e^{-it_0A_G}$ is uniform \cite{Chan20,GNR17}.

We index the entries of $A=[a_{kj}]\in\Cnn$ from $0$ to $n-1$. The set of  unitary $n\x n$ matrices is denoted by $\U(n)$, the set of  $n\x n$ matrices of determinant one is the \emph{special linear group} and is denoted by $\SL n$, and the set of invertible $n\x n$ matrices is denoted by $\gl n$.  Obviously $\U(n)$ and $\SL n$ are subgoups of $\gl n$, and we will sometimes have occasion to consider additional subgroups, such as the set of  real orthogonal $n\x n$ matrices, which is denoted by $\Q(n)$.    The \emph{max-norm} of $A$  is 
$\mnorm A =\underset{0\le k,j<n}{\max}\left|a_{kj}\right|$ and  the \emph{Frobenius norm} of $A$ is
$\Fnorm A = \displaystyle \sqrt{\tr \lp A^*A\rp} = \sqrt{\underset{0\le j,k < n} {\sum}\left|a_{kj}\right|^{2}}.$

\begin{defn}\label{d:uni-app}
A complex square matrix $A=[a_{kj}]$ is \emph{uniform} if there exists a nonnegative real number  $c$ such that $\left|a_{kj}\right|=c$ for all $k$ and $j$. A matrix $A\in\Cnn$ is \emph{unitarily apportionable} or \emph{$\U$-apportionable} if there exists a matrix $U\in \U(n)$ such that $UAU^*$ is uniform. In this case, $\mnorm{UAU^{*}}$ is called a \emph{unitary apportionment constant} and $U$ is called an \emph{apportioning matrix}. A matrix $A\in \Cnn$ is \emph{GL-apportionable} or \emph{generally apportionable} if  there exists a matrix  $M\in\gl n$ such that $MAM^{-1}$ is uniform.  In this case, $\mnorm{MAM^{-1}}$ is called a \emph{general apportionment constant} and $M$ is called an \emph{apportioning matrix}.  
\end{defn}

An apportionment constant is usually denoted by $\kappa$. If $A=[a_{kj}]\in\Cnn$ is uniform, then ${\mnorm A}= |a_{kj}|=\frac{\Fnorm A} n$. Since the Frobenius norm is unitarily invariant, the unitary apportionment constant is unique.
  
We have not defined \emph{specially apportionable}, because every generally apportionable matrix can be apportioned by a matrix in $\SL n$: If $MAM^{-1} = B$, then 
\[
(\det(M)^{-1/n} M)A(\det(M)^{1/n}M^{-1}) = B
\]
and $\det(M)^{-1/n} M\in\SL n$. Similarly, we may consider only special unitary matrices when studying $\U$-apportionment.

Example \ref{ex:app-not-u-app} shows that a matrix may be apportionable but not $\U$-apportionable.  Just as unitarily  apportionable matrices were defined to measure apportionability relative to $\U(n)$,  the apportionability of $A$ can be assessed relative to any subgroup of  $\gl n$.
\bigskip

We first study unitary apportionability and then consider general apportionability.  Determining whether a matrix is $\U$- or GL-apportionable can be challenging but every rank one matrix is $\U$-apportionable and we present an algorithm for finding a unitary matrix to apportion a rank one matrix in Section \ref{s:rank1}.  We also show there that a positive semidefinite matrix $H$ is $\U$-apportionable if and only if $\rank H\le 1$. In Section \ref{s:nec-translate} we present a condition on a matrix $A\in\Cnn$ that is sufficient to show that $A$ is not $\U$-apportionable; however, this condition is not necessary. In Section \ref{s:u-close} we define a function  that measures how far away from  $\U$-apportionable a matrix is and establish bounds on this function.
Connections with  Rosa's $\rho$-labelings of graphs, which include graceful labelings, are studied in Section \ref{s:graph-decomp}.  There we show that a loop-graph has a $\rho$-labeling if and only if a specific expansion of its adjacency matrix can be apportioned by a unitary matrix of a specific form.  In Section \ref{s:interlacing} we use gracefully labelled loop-graphs to generalize the well-known  Eigenvalue Interlacing  Inequalities.  General apportionment is studied in Section \ref{s:spectra}, where it is shown that most pairs of real numbers are not realizable as the spectrum of a $2\x 2$ apportionable matrix,  and Section \ref{s:M-kappa}, where the problem of finding an apportioning matrix is studied.  Section \ref{s:conclude} contains concluding remarks.  Some of the work in Section \ref{s:graph-decomp} relies on the Composition Lemma proved in \cite{G23}, and in Appendix \ref{appdx} we provide a proof of the Recovery Lemma, which is used in \cite{G23}.

The remainder of this introduction contains additional examples, terminology, and notation. 
 We use the notation $\ii$ for the imaginary unit and we work over the field of complex numbers unless otherwise indicated.  The $n\x n$ identity matrix is denoted by $I_n$ (or $I$ if $n$ is clear), the $n\x n$ all zeros matrix is denoted by $O_{n}$ (or $O$)  and 
the $n\x n$ all ones matrix is denoted by $J_n$ (or $J$). Let $E_{kj}$ be the $n\x n$ matrix with $(k,j)$-entry equal to 1 and all other entries 0 (note that $n$ must be specified when using $E_{kj}$).  Let $C_n$ denote  the circulant matrix with first row $[0,1,0,\dots,0]$, so $C_n$ is the adjacency matrix of a directed cycle spanning all $n$ vertices. 
Let $\bone=[1,\dots,1]^\top$, $\omega$ be a primitive $n$th root of unity,  $\bw=[1,\omega,\omega^2, \dots,\omega^{n-1}]^\top$, and $F_n=\frac 1 {\sqrt n}[\bw_0=\bone,\bw_1=\bw,\bw_2, \dots,\bw_{n-1}]$ where $\bw_j$ is the entrywise product of $j$ copies of $\bw$ for $j=2,\dots,n-1$. The matrix $F_n$ is called a \emph{discrete Fourier transform}  or \emph{DFT} matrix.  Since we index matrix entries by $0,\dots,n-1$, the $(k,j)$-entry of $F_n$ is $\frac{\omega^{kj}} {\sqrt n}$. DFT matrices are very useful for the study of apportionment and also provide a nice example of a family of uniform matrices for which the eigenvalues are completely known.

\begin{ex}\label{ex: DFT e-vals}
 It is known that every eigenvalue of  the $n\x n$ DFT matrix $F_n$  is one of $1,-1,\ii,-\ii$.  Furthermore the multiplicities are given in the next table. 
 \begin{table}[h!]\caption{Multiplicities of the eigenvalues $\pm 1, \pm \ii$ for $F_n$ \cite{McCP72}}\label{tab:DFT}
\begin{center} 
\noindent  \begin{tabular}{|c|c|c|c|c|}
\hline 

   $n$ &  $1$ & $-1$ &$-\ii$ & $\ii$ \\
 \hline 
$4k$ & 	$k+1$ &		$k$	& 	$k$	& $k-1$  \\
\hline 
$4k+1$ & 	$k+1$ &		$k$	& 	$k$	& $k$  \\
\hline $4k+2$ & $k+1$ &	$k+1$	& $k$	& $k$  \\
\hline $4k+3$ & $k+1$ &	$k+1$	& $k+1$	& $k$  \\
\hline 
\end{tabular}
\end{center}
\end{table}
\end{ex}

A  \emph{Hadamard matrix} is an $n\x n$  matrix  $H$ with every entry equal to  1 or $-1$ and such that $H^\top H = nI_n$. A Hadamard matrix of order $n$ is necessarily uniform and every eigenvalue has magnitude ${\sqrt n}$.  If there exists an $n\x n$ Hadamard matrix, then $n=1,2$ or $n\equiv 0\mod 4$. It is not known whether there exist Hadamard matrices of all orders of the form $n=4k$, but it is known that there is a Hadamard matrix for each $n=2^k$. If $H_m$ and $H_n$ are $m\x m$ and $n\x n$ Hadamard matrices, then $H_{mn}=H_m\otimes H_{n}$ is an $mn \x mn$ Hadamard matrix,  where $\otimes$ denotes the Kronecker product. Thus,  Hadamard matrices of order $n=2^k$  can be constructed as in the next example.    

\begin{ex}\label{ex-Hadamard}
Let $H_2=\mtx{1 & 1\\1 & -1}$ and 
define $H_{2^{k+1}}=H_2\otimes H_{2^k}$.
For a symmetric Hadamard matrix $H_n$, every eigenvalue is $\pm \sqrt n$.    Observe that $\tr(H_{2^k})=0$ and $H_{2^k}$ is symmetric, so $\spec(H_{2^k})=\lsb\lp-\sqrt {2^k}\rp^{(2^{k-1})},\lp\sqrt {2^k}\rp^{(2^{k-1})}\rsb$. By scaling $H_{2^k}$, any spectrum of the form $\lsb\lp -\lam\rp^{(2^{k-1})},\lam^{(2^{k-1})}\rsb$ can be realized by a $2^k\x 2^k$ uniform matrix. 
\end{ex}

There is a  close connection between the matrix apportionment problem and two classical graph decomposition problems. One of these problems (graceful labeling) is discussed in Section \ref{s:graph-decomp}.  The other problem consists in determining the existence of a decomposition of  $K_n$ (where $n$ is even) 
into $n-1$ overlapping copies of $K_{n/2,n/2}$ such that each edge of $K_n$ occurs as an edge in exactly $\frac{n}{2}$   distinct copies of the given $n-1$ copies of $K_{n/2,n/2}$.  This graph decomposition problem is well known to be equivalent to the problem of establishing the existence of a $n\times n$ Hadamard matrix. 
Its relation to apportioning follows from the observation that a symmetric $n\times n$ Hadamard matrix exists if and only if one of the $\frac{n}{2}$ diagonal matrices in the set
\[
\left\{ \diag(\underbrace{-1,\ldots,-1}_{k \text{ times}},\underbrace{1,\ldots,1}_{n-k \text{ times}}): 0 < k \leq \frac n2 \right\}
\]
 is {$\Q$}-apportionable. 
Thus determining whether or not one of the diagonal matrices above is {$\Q$}-apportionable must be at least as hard as establishing the existence of symmetric $n\times n$ Hadamard matrices.

 A \emph{complex Hadamard matrix} is an $n\times n$ matrix $H$ with complex entries of modulus 1 such that $HH^* = nI$.  An $n\x n$ complex Hadamard matrix is $\sqrt n \,U$ where $U\in\Cnn$ is both uniform and unitary.  The appropriate scalar multiple of the DFT matrix is a complex Hadamard matrix.
Thus, unlike real Hadamard matrices,  complex Hadamard matrices exist for  every order.

  We use $\spec(A)$ to denote the spectrum of $A$, i.e., the multiset of $n$ eigenvalues of $A$. The \emph{spectral radius} of $A$ is $\rho(A)=\max\{|\lam|:\lam\in\spec(A)\}$; $\rho$ can be used when $A$ is clear. The eigenvalues of a Hermitian matrix $H$ are real and denoted by $\lambda_{\max}(H)=\lambda_0(H)\ge \lambda_1(H)\ge \dots\ge \lambda_{n-1}(H)=\lambda_{\min}(H)$ (or $\lambda_{\max}=\lambda_0\ge \lambda_1\ge \dots\ge \lambda_{n-1}=\lambda_{\min}$ if $H$ is clear). 
 
Note that $\mnorm\cdot$  is a vector norm but not a matrix norm.  
We  use several matrix norms in addition to the Frobenius norm, many of which are defined in terms of the singular values of $A$, which we denote 
 by $\sigma_0(A)\ge\dots\ge \sigma_{n-1}(A)$.  If $A$ is a normal matrix, then $\sigma_0(A)=\rho(A)$ and $\{\sigma_0(A),\dots, \sigma_{n-1}(A)\}=\{|\lambda|: \lambda\in  \spec(A)\}$.  For a positive integer $p$  (or $p=\infty$), the Schatten-$p$ norm of $A\in\Cnn$  is
\[
\Snorm A p=\lp\sum_{0\le k<n} \sigma_k(A)^p\rp^{1/p};
\]
 $ \Snorm A 1=\sum_{0\le k<n} \sigma_k(A)$ is called the \emph{nuclear norm} and denoted by $\nnorm A$.  $ \Snorm A \infty= \sigma_0(A)$ is called the \emph{spectral norm} and denoted by $\tnorm A$ (because it is the matrix norm induced by the vector 2-norm). 
Since the Frobenius norm of $A$ is invariant under multiplication of $A$ by a unitary matrix,
$\Fnorm A=\Snorm A 2;$ note that the Frobenius norm of $A$ is the vector 2-norm of $A$ viewed as an $n^2$-vector. 

\section{Rank one matrices are $\U$-apportionable}\label{s:rank1} 

 The only rank zero matrix in $\Cnn$ is $O_n$, and it is uniform.  In this section we show that every rank one matrix is $\U$-apportionable  and  provide an algorithm for finding a unitary apportioning matrix and similar uniform matrix.  
The situation changes dramatically when the  rank is of two or more.   In this section we also show that positive semidefinite matrices of rank two or more are not $\U$-apportionable.  Section \ref{s:spectra} provides examples of spectra that cannot be realized by apportionable matrices.

Lemma \ref{lem: r1 canonical} and Theorem \ref{t:u-apport-exist} establish that Algorithm \ref{alg1} produces the claimed results; to assist in making  connections between the proofs and the algorithm,   we identify algorithm steps by number within the proofs of Lemma \ref{lem: r1 canonical} and Theorem \ref{t:u-apport-exist}.

Recall that for any vectors $\bv,\bw\in\Cn$ with $\tnorm{\bv}=\tnorm{\bw}$, there exists a unitary matrix $U$ such that $U\bv=\bw$. This can be accomplished, for example, by the \emph{Householder matrix  defined by $v$ and $w$,} $I_n-\frac{\bu \bu^*}{\bu^*\bv}$, where $\bu=\bv-\bw$ (cf. \cite{CplxHouse}).   The (complex) sign function is  $\sgn(z)=\frac z{|z|}$ for all nonzero $z\in\C$ and $\sgn(0)=1$.  Let $\be_i$ denote the $i$-th standard basis vector.

\begin{lem}\label{lem: r1 canonical}
Let $n\geq 2$ and $A\in\C^{n\times n}$. If $\rank A = 1$, then $A$ is unitarily similar to $\gamma\be_0(\alpha\be_0^\top+ \beta\be_1^\top)$, where $|\gamma| = 1$ and $\alpha,\beta\in\R$ are nonnegative.
\end{lem}
\begin{proof}
Suppose that $\rank A = 1$.  This implies that $\spec(A)=\{\lambda, 0^{(n-1)}\}$ where $\lam=\tr A$ (Step 1 in Algorithm \ref{alg1}); $\lam$ may be nonzero or zero.  Furthermore, $A = \bx \by^*$ for some $\bx,\by\in \C^n$ with $\bx,\by\ne\bzero$ and $\|\bx\|_2 = 1$  (Step 2).  Let $H_1$ be a {Householder} matrix such that $H_1\bx = \be_0$ and let $\bz=H_1\by$  (Steps 3 and 4). This implies $(\bz)_0=\ol\lam$.
Define $\hat\bz=\bz-{\ol\lam}\be_0$, so $(\hat\bz)_0=0$. If $\hat\bz=\bzero$, then let $ H_2=I_n$.  Otherwise, let $H_2$ be  {the Householder matrix defined by $\hat\bz$ and $\sgn({\ol\lam})\tnorm{\hat\bz}\be_1$  (Step 5).
Thus $H_2\hat\bz = \sgn({\ol\lam})\tnorm{\hat\bz}\be_1$,  $H_2\be_{0}=\be_{0}$, and $H_2 H_1 \by = \sgn({\ol\lam})(|\lam| \be_0 + \tnorm{\hat\bz}\be_1)$.  Hence, 
 $H_2 H_1 A H_1^* H_2^* = {\sgn(\lam)}\be_0(|\lam| \be_0 + \tnorm{\hat\bz}\be_1)^\top$ has the required form.}
\end{proof}

\begin{algorithm}[!h]
\begin{algorithmic}[1]
\caption{\label{alg1}  Given $A\in\Cnn$ with $\rank A=1$, construct uniform matrix $B$  and unitary matrix $V$ such that $B=VAV^*$.}

\STATE   $\lambda=\tr A$.
\STATE   Factor $A$ as $A=\bx\by^*$ where  $\tnorm{\bx}=1$.
\STATE  $H_1=I_n-\frac{\bu\bu^*}{\bu^*\bx}$ where $\bu=\bx-\be_0$ (assuming $\bx\ne \be_0$, else $H_1=I_n$).
\STATE $\bz=H_1\by$.
\STATE   $H_2=I_n-\frac{\hat\bu\hat \bu^*}{\hat\bu^*\hat \bz}$ where  $\hat \bz=[0,z_1,\dots,z_{n-1}]^\top$ and ${\hat\bu}=\hat \bz-\sgn(z_0)\tnorm {\hat\bz}\,\be_1$  \\ \hspace{7pt}(assuming $\hat \bz \ne\bzero$ and $\hat \bz\ne\sgn(z_0) \tnorm {\hat\bz}\,\be_1$, else $H_2=I_n$).

\IF{$\lambda= 0$}
\STATE   $V=F_nH_2H_1$.
\STATE  $(B)_{kj}=VAV^*=\frac{\Fnorm A \, \omega^{-j}}n$.
\ENDIF

\IF{$\lambda\ne 0$ {\bf and} $\Fnorm A =|\tr A|$}
    \STATE   $V=F_nH_2H_1$.
        \STATE   $B=VAV^*=\frac{\Fnorm A}n J_n$.
\ENDIF

\IF{$\lambda\ne 0$ {\bf and} $\Fnorm A \ne |\tr A|$}
\STATE   $\bu_0=\frac 1 {\sqrt n} \bone$.
\IF{$n$ is even}
    \STATE   $\bu_1=\frac 1 {\sqrt n}[\ii,-\ii,\dots,\ii,-\ii]^\top$.
 \ELSE
    \STATE   $\bu_1=  [(1-n)a, a + b\ii, a - b\ii, a + b\ii,\ldots, a - b\ii]^\top$ where\vspace{-8pt} \[
a = \frac{1-\sqrt{\beta^2+1}}{(n-1) \sqrt{n} \beta}
\qquad\mbox{ and }\qquad b 
= \sqrt{(n-1)^{-1} - n a^2}
.\vspace{-8pt}\] 
 \ENDIF
 \STATE   Construct a unitary matrix $U=[\bu_0, \bu_1, \bu_2,\dots,\bu_{n-1}]$ (e.g., by extending $\{\bu_0,\bu_1\}$ to a basis and\\ \hspace{7pt}  applying the standard Gram-Schmidt process). 
 \STATE   $V=U H_2H_1$
 \STATE   $B=VAV^*$
\ENDIF

\end{algorithmic}
\end{algorithm}

  {Note that this algorithm is intended to summarize the constructive method of proof.  For accurate implementation in decimal arithmetic, it is important to apply well-known numerical techniques to minimize errors.}
 
\begin{thm}\label{t:u-apport-exist}
Let $A\in \C^{n\times n}$. If $\rank A = 1$, then $A$ is $\U$-apportionable.
\end{thm}
\begin{proof}
The claim clearly holds for $n = 1$, so suppose $n\geq 2$. Assume that $\rank A = 1$. By Lemma \ref{lem: r1 canonical} we may assume that $A$ is of the form $\gamma\be_0(\alpha\be_0^\top + \beta\be_1^\top)$, where $|\gamma| = 1$ and $\alpha,\beta\in\R$ are nonnegative.  Observe that if $\alpha = 0$ or $\beta =0$, then $F_n A F_n^*$ is uniform (Steps 6-13 of Algorithm \ref{alg1}). So suppose that $\alpha$ and $\beta$ are positive real numbers. Since $\U$-apportionability is invariant under scaling we may assume that $\gamma = 1$ and $\alpha = 1$. 

 Let $U$ be a unitary matrix whose first two columns are $\bu_0$ and $\bu_1$. Then 
\[
UAU^* = ( U\be_0\, (\be_0^\top + \beta\be_1^\top)U^*) = \bu_0\, (\bu_0 + \beta\bu_1)^*.
\]
Thus $UAU^*$ is uniform if and only if $\bu_0$ and $\bu_0 + \beta\bu_1$ are uniform.
If $n$ is even and  $U$ is a unitary matrix whose first two columns are $\bu_0=\frac 1 {\sqrt n} \bone$ and $\bu_1=\frac 1 {\sqrt n}[\ii,-\ii,\ii,\ldots, -\ii]^\top$, then $\bu_0$ and $\bu_0 + \beta\bu_1$ are uniform (Steps 14-17).

Now consider the case where $n$ is odd. Let $U$ be a unitary matrix whose first two columns are $\bu_0=\frac 1 {\sqrt n} \bone$ and
\[
\bu_1=  [(1-n)a, a + b\ii, a - b\ii, a + b\ii,\ldots, a - b\ii]^\top,
\]
where $b = \sqrt{(n-1)^{-1} - na^2}$ and $a = \frac{1 - \sqrt{\beta^2 + 1}}{(n-1)\sqrt{n}{ \beta}}$ (Steps 14-15 and 18-19). Then it is not difficult, albeit rather tedious, to verify algebraically that $\bu_0 + \beta\bu_1$ is uniform, i.e.,
\[
\lab\frac1{\sqrt{n}} + \beta(1-n)a\rab= \sqrt{\frac{\beta^2+1}n}= \lab\frac1{\sqrt{n}} + \beta\lp a + \sqrt{ (n-1)^{-1} - na^2}\ii\rp\rab = \lab\frac1{\sqrt{n}} + \beta(a + b\ii)\rab. \qedhere
\]
\end{proof}

Recall that a matrix $H$ is \emph{positive semidefinite (PSD)} if  and only if $H$ is Hermitian and   $\lambda\ge 0$ for every eigenvalue $\lambda$ of $H$. 
An $n\times n$ matrix of the form $H = X^*X$ for some $X\in \C^{m\x n}$ is called a \emph{Gram matrix}.    It is well known that a matrix is a Gram matrix if and only if it is positive semidefinite.  In fact,  the least $d$ such that a PSD matrix $H$ can be expressed as $X^*X$ with $X\in\C^{d\x n}$ is the rank of $H$.

\begin{prop}\label{p:PSD-conj}\label{p:PSD}
Let $H$ be a PSD matrix. If $C$ is nonsingular and $C^*HC$ is uniform, then $\rank H\le 1$.  
Furthermore, $H$ is  $\U$-apportionable if and only if  $\rank H\leq 1$.
\end{prop}
\begin{proof}
If $H=O$ then $\rank H=0$.  So assume $H\ne O$.  Let $C$ be nonsingular and suppose $B=C^*HC$ is uniform. Since $H$ is PSD, $B$ is also PSD and there exists a matrix $R$ such that $B=R^*R$.  
Let $\br_k$ denote  column $k$ of $R$. Since $ B$ is uniform, $\tnorm{\br_k}=\tnorm{\br_j}$ and $|\br_k^*\br_j| = \tnorm{\br_k}\tnorm{\br_j}$ for any row indices $k$ and $j$. Thus equality holds for the Cauchy-Schwarz inequality applied to any pair of  columns of $R$. This implies $\rank R = 1$, and thus $\rank H= 1$.

If $\rank H\le 1$, then  
$H$ is $\U$-apportionable by Theorem \ref{t:u-apport-exist}.
 Now suppose  that $H$ is $\U$-apportionable, so there exists a unitary matrix $U$ such that $U^*HU$ is uniform.  Then  $\rank H\le 1$. 
\end{proof}

The next example shows that a PSD matrix $H$ with $\rank H=2$ may be GL-apportionable, demonstrating the existence of a matrix that is apportionable but not $\U$-apportionable. 
\begin{ex}\label{ex:app-not-u-app}
 Observe that
\[
A=\left[\begin{array}{ccc}
1 & 1 & -1\\
1 & 1 & 1\\
-1 & -1 & 1
\end{array}\right]
=
\left[\begin{array}{ccc}
1 & 1 & 1\\
-1 & -1 & 0\\
0 & -1 & -1
\end{array}\right]
\left[\begin{array}{ccc}
0 & 0 & 0\\
0 & 1 & 0\\
0 & 0 & 2
\end{array}\right]
\left[\begin{array}{ccc}
1 & 0 & 1\\
-1 & -1 & -1\\
1 & 1 & 0
\end{array}\right]
\]
is a uniform matrix with distinct nonnegative  eigenvalues, $H=\mtx{0 & 0 & 0\\
0 & 1 & 0\\
0 & 0 & 2}$ is PSD, and $\rank H= 2$. \end{ex}

It is interesting to note that while the previous example shows that $\{1,2,0\}$ is realizable as the spectrum of a uniform matrix, Theorem \ref{2x2-real-spec} implies $\{1,2\}$ is not realizable. 
 
\section{Necessary condition for $\U$-apportionability}
\label{s:nec-translate}
 In this section we establish a necessary condition on $|c|$ for a translation $A+cI$ to be $\U$-apportionable. This condition   is used to show  that for a given $n\ge 2$, a positive fraction of  $n\x n$ matrices are not  $\U$-apportionable.  

\begin{thm}\label{t:necessary}
Let $n\ge 2$, $A=[a_{kj}]\in\Cnn$, and $c\in\C$. Suppose $A+cI_n$ is $\U$-apportionable.  Then 
\[ |c|\le \frac{ \sum_{k=0}^{n-1} |a_{kk}|} n+\sqrt{\lp\frac{ \sum_{k=0}^{n-1} |a_{kk}|}n\rp^2+\frac{\Fnorm A ^2}{n(n-1)}}\]
\end{thm}
\bpf Let $B=A-A\circ I_n$.  Since $A+cI_n$ is $\U$-apportionable, the $\U$-apportionment constant for $A+cI_n$ is 
\[\kappa=\frac{\Fnorm{A+cI}}n=\frac{\Fnorm{A\circ I_n+cI+B}}n=\frac 1 n \sqrt{\sum_{k=0}^{n-1} |a_{kk}+c|^2+\Fnorm{B}^2}. \]
 Since $A + cI$ is $\U$-apportionable, there exists a unitary matrix $U$ such that $U(A+cI)U^*$ is uniform. Observe that $U(A+cI)U^*=UAU^*+cI$.  Thus $|\lp UAU^*\rp_{kj}|=\kappa$ for all $k\ne j$. So
\[ \Fnorm A ^2 = \Fnorm{UAU^*}^2\ge \sum_{k=0}^{n-1}\sum_{j\ne k}|\lp UAU^*\rp_{kj}|^2=n(n-1)\kappa^2.\]
Replacing $\Fnorm A^2$ and $\kappa^2$ and rearranging yields the following.
\bea
\sum_{k=0}^{n-1}|a_{kk}|^2+\Fnorm B^2&\ge &\frac{n-1}n\lp  \sum_{k=0}^{n-1} |a_{kk}+c|^2+\Fnorm{B}^2\rp.\\
\sum_{k=0}^{n-1}|a_{kk}|^2 +\frac 1 n\Fnorm B^2&\ge &\frac{n-1}n  \sum_{k=0}^{n-1} (|a_{kk}|-|c|)^2\\
&= &\frac{n-1}n \lp \sum_{k=0}^{n-1} |a_{kk}|^2 + \sum_{k=0}^{n-1}|c|^2- 2\sum_{k=0}^{n-1} |a_{kk}||c| \rp.\\
\frac 1 n\sum_{k=0}^{n-1}|a_{kk}|^2 +\frac 1 n\Fnorm B^2&\ge & (n-1) |c|^2- \lp\frac{2(n-1)}n \sum_{k=0}^{n-1} |a_{kk}|\rp|c|.\\
0&\ge & (n-1) |c|^2- \lp\frac{2(n-1)}n \sum_{k=0}^{n-1} |a_{kk}|\rp|c|-\frac 1 n \Fnorm A.\eea
Since $n\ge 2$, we can apply the quadratic formula to obtain
\[ |c|\le \frac{ \sum_{k=0}^{n-1} |a_{kk}|} n+\sqrt{\lp\frac{ \sum_{k=0}^{n-1} |a_{kk}|}n\rp^2+\frac{\Fnorm A ^2}{n(n-1)}}.\qedhere\]
\epf

 Note that it follows from Theorem  \ref{constant-spec} below that translation also does not preserve GL-apportionability.

  An equivalent form of Theorem \ref{t:necessary}, 
\beq\label{eq:c-test}\exists c, |c|>  \frac{\sum_{k=0}^{n-1}|a_{kk}-c|}n+\sqrt{\lp\frac{\sum_{k=0}^{n-1}|a_{kk}-c|}n\rp^2+\frac{\Fnorm {A-cI}^2}{n(n-1)}}  \mbox{ implies $A$ is not $\U$-apportionable},\eeq
 gives a condition that can be used to show a matrix is not $\U$-apportionable.  We apply this
 to show that 
it is not the case that almost all $n\x n$ matrices are $\U$-apportionable. 
Consider the complex vector space $\Cnn$ of $n\x n$ matrices with max-norm (or equivalently, $\C^{n^2}$ with $\infty$-norm  where a matrix  represented by the vector of its  entries). 
Define the region 
$T_n=\{A\in \Cnn: \mnorm{A-\frac 3 4 I_n} \le \frac 1 4
\}$.

\begin{prop}\label{c:notU}
For $n\ge 2$, no matrix in $T_n$ is $\U$-apportionable.   
\end{prop}
\bpf For $c=\frac 3 4$,
\bea  \frac{\sum_{k=0}^{n-1}|a_{kk}-c|}n+\sqrt{\lp\frac{\sum_{k=0}^{n-1}|a_{kk}-c|}n\rp^2+\frac{\Fnorm {A-cI}^2}{n(n-1)}} &\le& \frac 1 4+\sqrt{\lp\frac 1 4\rp^2+\frac{n\lp\frac 1 4\rp^2+n(n-1)\lp\frac 1 4\rp^2}{n(n-1)}}\\
&\le& \frac 1 4+ \sqrt{\frac 1{16}+\frac 1 {16(n-1)}+\frac 1{16}}<\frac 3 4.\eea
 The result now follows by applying condition \eqref{eq:c-test}. \epf

 \begin{rem}
    Since $T_n$ has nonzero Lebesgue measure and is scalable, it is not the case that almost all $n\x n$ matrices are $\U$-apportionable.
\end{rem}

Although \eqref{eq:c-test} is a sufficient condition, it is not a necessary condition, as we can see by examining a connection between $\U$-apportionabilty and equiangular lines.  The unit vectors $\bx_0,\ldots,\bx_{d-1}\in\C^n$ are \emph{equiangular with angle $\theta$} if $|\bx_i^* \bx_j| = \cos(\theta)$  for all $0\le i< j\le d-1$; note that ${\bx_j}^*\, \bx_j = 1$ for $j=0,\dots,d-1$. 
 Given $B \in \mathbb{C}^{d \times n}$  with unit length columns, consider determining the existence of a unitary matrix $U$ whose   right action on $B$ results in (unit length) equiangular  columns (with angle $\theta$).
Note that for any $c\in\C$, $U^*(B^*B+cI_n)U=(BU)^*(BU)+cI_n$.  Thus $BU$ has equiangluar columns with angle $\theta$ if and only if $
    |\lp(BU)^*(BU)\rp_{ij}|= \cos(\theta)\mbox{ and }\lp(BU)^*(BU)\rp_{jj}=1$ if and only if 
     \[U^*\lp B^{*}B-I_{n}+
 \cos(\theta) I_n 
 \rp U\ \mbox{ is uniform}\] 
So the matrix $ B^*B+
\lp \cos(\theta) -1\rp I_n$ 
 is $\U$-apportionable  if and only if there exists a unitary matrix $U$ such that the columns of $BU$ are equiangular  with angle $\theta$.

Observe that $B^*B-I$ has zeros on the diagonal and   $\Fnorm{B^*B-I}=\Fnorm{U(B^*B-I)U}$ for any unitary matrix $U$.  So if $B^*B-I$ is $\U$-apportionable, the $\U$-apportionment constant is $\kappa=\sqrt{\frac{\Fnorm{B^*B-I}^2}{n(n-1)}}$ and we are trying to apportion $A=B^*B+(\kappa-1)I$.  However, we cannot use Theorem \ref{t:necessary} to show this is not $\U$-apportionable, because we see that for any $c\in \C$:
\bea \frac{\sum_{k=0}^{n-1}|a_{kk}-c|}n+\sqrt{\lp\frac{\sum_{k=0}^{n-1}|a_{kk}-c|}n\rp^2+\frac{\Fnorm {A-cI}^2}{n(n-1)}} &=& |\kappa-c|+\sqrt{|\kappa-c|^2+\frac{\Fnorm {A-cI}^2}{n(n-1)}}\\
&\ge&|\kappa-c|+\sqrt{|\kappa-c|^2+\kappa^2}\\
&\ge& (|c|-\kappa) +\kappa =|c|
\eea

It is well known that there is a limit to the number of equiangluar lines in a given dimension $d$.  In the case $d=2$, no set of $4$ lines is  equiangular.  Thus for any $2\x 4$ matrix $B$, $A=B^*B+(\sqrt{\frac{\Fnorm{B^*B-I}^2}{n(n-1)}}-1)I$ is not $\U$-apportionable, but no $c$ satisfies the hypothesis of  \eqref{eq:c-test}.

 \begin{quest}\label{q:certificate} Can we efficiently certify that a matrix is not $\U$-apportionable? 
\end{quest}

 \section{Measure of closeness to uniformity}\label{s:u-close}
 In this section we  define a function  that measures how far away from  $\U$-apportionable a matrix is and establish bounds on this function.
\begin{defn}
For $A\in\Cnn$, define
\[
u(A)=\inf_{U\in \U(n)}\mnorm{UAU^*},
\]
The \emph{unitary apportionment gap} of an arbitrary matrix $A\in \mathbb{C}^{n\times n}$ is $\left|u(A)-\frac{\left\Vert A\right\Vert _{\text{F}}}{n}\right|$
\end{defn}
Observe that $A$ is $\U$-apportionable if and only if $u(A)=\frac{\Fnorm A}n$ and the unitary apportionment gap is zero.
 Note that in defining the unitary apportionment gap, we are using the fact that the unitary apportionment constant is unique.
\begin{rem}\label{r:U-compact}
    Since the set of $n\x n$ unitary matrices is compact, the infimum is actually the minimum:
    \[
u(A)=\min_{U\in \U(n)}\mnorm{UAU^*}.
\]
\end{rem}
Next we establish bounds on $\mnorm{UAU^*}$ and thus on $u(A)$.  Recall that the nuclear norm of $A$ is $\nnorm A=\sum_{k=1}^n \sigma_k(A)$ and the spectral norm of $A$ is $\tnorm A=\sigma_0(A)$ (where $\sigma_k(A)$ is the $k$th singular value of $A$).
\begin{prop}\label{p:u-bds2} Let $A=[a_{ij}]\in\Cnn$.  Then 
\[\frac{\Fnorm A}n\le \mnorm{UAU^*}  \le \tnorm A\] 
The lower bound is realized by a uniform matrix unitarily similar to $A$ if $A$ is $\U$-apportionable.  The upper bound is realized by a diagonal matrix unitarily similar to $A$ when $A$ is normal.
\end{prop} 
\bpf 
 It is known that for any matrix $A\in \Cnn$, $\mnorm A\le \tnorm A$ (e.g., this is immediate from \cite[Theorem 5.6.2(d)]{HJ}).
Thus  $\mnorm{UAU^*}  \le \tnorm{UAU^*}=\tnorm A$. 
 For any matrix $B\in \Cnn$, $\Fnorm B\le \sqrt{n^2\lp\max_{k,j}|b_{kj}|\rp^2}=n\mnorm B$ with equality if and only if $B$ is uniform. 
Thus  $\Fnorm A\le n \mnorm{UAU^*}$ since  $\Fnorm{UAU^*}=\Fnorm A$ for any unitary matrix $U$.  
The last two statements are now immediate. 
\epf

\begin{prop}\label{p:u-bds} Let $A=[a_{ij}]\in\Cnn$.  Then 
\beq\label{eq:u-bds} \frac{\Fnorm A}n\le u(A)  \le \tnorm A\eeq
and $\frac{\Fnorm A}n= u(A)$  if and only if $A$ is $\U$-apportionable. 
If $A$ is normal, then $u(A) \le \frac{\nnorm A}n$.
All bounds are sharp.
\end{prop}

\bpf Equation \eqref{eq:u-bds} is immediate from Proposition \ref{p:u-bds2},  where it is also shown that equality in the first inequality requires $A$ to be $\U$-apportionable.

Now assume $A$ is normal and recall $F_n$ is the DFT. Since $A$ is normal, there exists a unitary matrix $V$ such that $V^*AV=\diag(\lam_1,\dots,\lam_n)$ where $\spec(A)=\{\lam_1,\dots,\lam_n\}$.  The singular values of $A$ are the absolute values of the eigenvalues of $A$ because $A$ is normal, so
\[  \big |({F_nV^*AVF_n^*})_{\ell j} \big | =\Big|\sum_k \frac {\omega^{\ell k}}{\sqrt n}\lam_k \frac {\ol{\omega^{kj}}}{\sqrt n}\Big|\\
\le\frac 1 n \sum_k |\omega^{\ell k}| |\lam_k| |\ol{\omega^{jk}}|\\
=\frac 1 n\sum_k  \sigma_k \\
=\frac {\nnorm A}n.\]
  Thus $u(A)\le \frac{\nnorm A}n$.
The upper  bounds are equality  when $A$ is the identity matrix, which is normal. \epf

\begin{quest}
    Is  upper bound on $u(A)\le \tnorm A$ sharp for any $A$ with $\frac{\nnorm A}n <\tnorm A$?  
\end{quest} 

\section{Graph labelings and $\U$-apportionment}\label{s:graph-decomp} 

One of the inspirations for writing this paper is a connection between apportionment and graph labelings introduced by Rosa in \cite{R1967}. Rosa demonstrated that some of these labelings produce cyclic decompositions of the complete graph. We show that these cyclic decompositions can naturally be described as an apportionment problem using a reasonable choice of unitary matrix and adjacency matrices.  

We first provide the necessary background information. The term \emph{graph} will always refer to a simple graph $G$, i.e., $G$ does not have any loops or multi-edges. A \emph{loop-graph} $\mathfrak G$ is a graph that allows loops but not multi-edges. We will always label the vertices of a graph on $n$ vertices with the set $\{0,\ldots,n-1\}$. The adjacency matrix of a graph $G$ or loop-graph $\mathfrak{G}$ is denoted $A_G$ or $A_{\mathfrak{G}}$, respectively.  The edge set of a graph $G$ or loop-graph $\mathfrak G$ is denoted by $E_G$ or $E_{\mathfrak G}$.

We denote the complete graph on $n$ vertices by $K_n$. The complete loop-graph $\mathfrak{K}_n$ is the loop-graph obtained from $K_n$ by adding a loop to each vertex. Define $\phi:V(K_n)\to V(K_n)$ to be the graph isomorphism that maps $i \mapsto i + 1 \mod n$. Note that $\phi$ is also a graph isomorphism of $\mathfrak{K}_n$.   Recall that $C_n$  denotes the  $n\times n$ (cyclic) permutation matrix corresponding to $\phi$, i.e., $C_n \be_i = \be_{i+1\mod n}$. 
The \emph{length} of an edge $\{i,j\}$ in $\mathfrak{K}_n$ is $\min\{|i-j|, n- |i-j|\}$. Observe that $\phi$ does not change the length of an edge and when  $n$ is odd, $\mathfrak{K}_n$ consists of 
$n$ edges of length $i$ for $i = 0,\ldots, \frac{n-1}{2}$.

Given a loop-graph $\mathfrak{G}$, a $\mathfrak{G}$-\emph{decomposition} of $\mathfrak{K}_n$ is a set $\Delta = \{\mathfrak{G}_1,\ldots,\mathfrak{G}_t\}$ of subgraphs of $\mathfrak{K}_n$ such that each $\mathfrak{G}_i$ is isomorphic to $\mathfrak{G}$ and the edge sets of the loop-graphs in $\Delta$ partition the edge set of $\mathfrak{K}_n$. The $\mathfrak{G}$-decomposition $\Delta$ is \emph{cyclic} if $\phi(\mathfrak{G}_i) = \mathfrak{G}_{(i+1)\mod n}$ for each $\mathfrak{G}_i$.  Let $\mathfrak{G}$ be a loop-graph with $m$ edges. An injective function $f: V(\mathfrak{G}) \to \{0,\ldots, 2m-2\}$ is a \emph{$\rho$-labeling} of $\mathfrak{G}$ provided 
\[
\big\{\min\{|f(u) - f(v)|,\, 2m+1-|f(u) - f(v)|\}: \{u,v\}\in E(\mathfrak{G})\big\} = \{0,\ldots, m-1\}.
\]
Thus a $\rho$-labeling is an embedding from $\mathfrak{G}$ to $\mathfrak{K}_ 2m-1$ such that the image of $\mathfrak{G}$ has exactly one edge  of length $i$ for $i = 0,\ldots, m-1$ (an \emph{embedding} of $\mathfrak G$ in $\mathfrak G'$ is an injective mapping of vertices that maps edges of $\mathfrak G$ to edges of $\mathfrak G'$).

The definitions in the preceding two paragraphs were originally given for simple graphs. We have modified the original definitions to fit the setting of loop-graphs. This small change simplifies the connection to apportionment. The following theorem is a loop-graph version of Theorem 7 in \cite{R1967}.

\begin{thm}\label{Thm: Rosa Rho Labeling}{\rm\cite{R1967}}
Let $\mathfrak{G}$ be a loop-graph with $m$ edges, consisting of  $m-1$ non-loop edges and one loop. Then there exists a cyclic $\mathfrak{G}$-decomposition of the complete loop-graph $\mathfrak{K}_{2m - 1}$ if and only if $\mathfrak{G}$ has a $\rho$-labeling.
\end{thm}

Let $\mathfrak{G}$ be a loop-graph with $n$ vertices, $n-1$ non-loop edges and one loop edge. Let $f: V(\mathfrak{G}) \to V(\mathfrak{K}_{2n-1})$ be an 
embedding and define $\widehat{\mathfrak G}$ to be the  
loop-graph with vertex set $\{0,\dots,2n-1\}$ and edge set $f(E_{\mathfrak G})$. In the special case that  $f(i)=i$ for all $i$, 
 the $(2n-1)\times (2n-1)$ adjacency matrix of $\widehat{\mathfrak G}$  would be $A_{\mathfrak{G}}\oplus O$.  
In the more general setting, a permutation similarity corresponds to a relabeling of vertices, so  there exists a permutation matrix $P$ such that the adjacency matrix of $\widehat{\mathfrak G}$ is  $A = P(A_{\mathfrak{G}}\oplus O)P^T$. Observe that there exists a cyclic $\mathfrak{G}$-decomposition of $\mathfrak{K}_{2n-1}$ if and only if
\begin{equation}\label{Uniformity}
\sum_{0\le i<2n-1}C_{2n-1}^{i}AC_{2n-1}^{-i}=J_{2n-1}.
\end{equation}
This observation can be thought of as a matrix version of Theorem \ref{Thm: Rosa Rho Labeling}. With some minor modifications, we can directly connect these ideas to apportionment.

Let $\omega$ denote a primitive $(2n-1)$-th root of unity and let $\bw= [1,\omega, \omega^2,\dots,\omega^{2n-2}]^\top$ (for convenience we may assume $\bw$ is  the second column of  $\sqrt{2n-1}\,F_{2n-1}$). Define $U_n\in \mathbb{C}^{(2n-1)^{2}\times(2n-1)^{2}}$
to be the $(2n-1)\times(2n-1)$ block matrix
whose  $(i,j)$-block is the  $\left(2n-1\right)\times\left(2n-1\right)$ matrix

\begin{equation}\label{Unitary_Matrix}
U_n(i,j) = \frac{C_{2n-1}^{j} \diag(\bw)^i}{\sqrt{2n-1}}, \quad\text{ for } 0\le i,\,j<2n-1.
\end{equation}
Observe that the $(i,k)$-block of $U_n U_n^*$ is 
\[
  \lp U_nU_n^*\rp(i,k) = \sum_{j=0}^{2n-2} U_n(i,j) U_n(k,j)^* = \frac 1 {2n-1} \sum_{j=0}^{2n-2} C_{2n-1}^{j} \diag \lp 1,\omega^{i-k},\omega^{2(i-k)},\dots,\omega^{(2n-2)(i-k)}\rp C_{2n-1}^{-j}.
\]
 If $i=k$, then $\diag\lp 1,\omega^{i-k},\omega^{2(i-k)},\dots,\omega^{(2n-2)(i-k)}\rp=I_{2n-1}$, and so $\lp U_nU_n^*\rp(i,i) = I_{2n-1}$. If $i\not=k$, then
\[
\lp U_nU_n^*\rp(i,k) = \sum_{j=0}^{2n-2} \diag\lp \omega^{j(i-k)},\omega^{(j+1)(i-k)},\dots,\omega^{(2n-2)(i-k)},1, \dots \omega^{(j-1))(i-k)}\rp = O_{2n-1}.
\]
Thus $U_n$ is unitary. Using $U_n$ we define the following  matrix representation of the symmetric group  $S_{2n-1}$ on $2n-1$ elements:   
\[
\mathfrak{U}_{n}:= \left\{ U_{n} (I_{2n-1} \otimes P) U_n^* : P \text{ is a $(2n-1) \times (2n-1)$ permutation matrix}\right\}.
\]
Note that $\mathfrak{U}_n$ is a subgroup of the unitary group  $\U((2n-1)^2)$. Define the \emph{cyclic blowup matrix} of $\mathfrak{G}$ to be the $(2n-1)^2 \times (2n-1)^2$ Hermitian matrix
\[
H_{\mathfrak{G}} = U_n(I_{2n-1} \otimes (A_{\mathfrak{G}} \oplus O_{n-1}))U_n^*.
\]
With these definitions in hand, we are now ready for the main result of this section.

\begin{thm}\label{t:Con2apport}
Let $\mathfrak{G}$ be a loop-graph with $n-1$ edges and one loop. Then $\mathfrak{G}$ has a $\rho$-labeling if and only if $H_{\mathfrak{G}}$ is $\mathfrak{U}_n$-apportionable. In this case, the $\U$-apportionment constant of $H_{\mathfrak{G}}$ is $\kappa = (2n - 1)^{-1}$.
\end{thm}
\begin{proof}
Assume that $\mathfrak{G}$ has a $\rho$-labeling $f: V(\mathfrak{G})\to V(\mathfrak{K}_{2n-1})$. Let $P$ be the permutation matrix such that the $(2n-1)\times(2n-1)$ adjacency matrix of the image of $\mathfrak{G}$ under $f$ is $A = P(A_{\mathfrak{G}}\oplus O)P^*$. Let $Q = U_n(I_{2n-1}\otimes P)U_n^*$. Note that $Q\in\mathfrak{U}_n$. Direct calculation gives
\[
QH_{\mathfrak{G}} Q^* = U_n(I_{2n-1}\otimes A)U_n^*
\]
and so the $(i,k)$-block of $QH_{\mathfrak{G}} Q^*$ is 
\[
\sum_{j=0}^{2n-2} U(i,j)A U(k,j)^* = \frac1{2n-1} \sum_{j=0}^{2n-2} C_{2n-1}^j \diag(\bw)^i A \diag(\bw)^{-k} C_{2n-1}^{-j}.
\]
Since $\mathfrak{G}$ has $n-1$ non-loop edges and one loop, Theorem \ref{Thm: Rosa Rho Labeling} implies that each term in this sum has disjoint support if and only if $f$ is a $\rho$-labeling. This holds if and only if the $(i,k)$-block of $QH_{\mathfrak{G}} Q^*$ is uniform with apportionment constant $(2n-1)^{-1}$  since the zero-nonzero pattern of each term in this sum matches that of the corresponding term in (\ref{Uniformity}). 
\end{proof}

Rosa also introduced other labelings in \cite{R1967}, in particular graceful labelings. Similar versions of Theorem \ref{t:Con2apport} can be made for each these labelings since they all imply the existence of $\rho$-labelings.
A graceful labeling of an $n$-vertex graph is a $\rho$-labeling where the previously mentioned injective function $f: V(\mathfrak{G}) \to \{0,\ldots, 2n-2\}$ is such that its  image is $\mathbb{Z}_n:=\{0,\ldots, n-1\}$. In a graceful labeling, the subset of vertices being acted upon by the permutation group is reduced from $2n-1$ to $n$, thereby further reducing the apportioning unitary subgroup from a representation of $S_{2n-1}$ as
\[
\mathfrak{U}_{n}:= \left\{ U_{n} (I_{2n-1} \otimes P) U_n^* : P \text{ is a $(2n-1) \times (2n-1)$ permutation matrix}\right\}
\]
to a representation of $S_n$ as
\[
\mathfrak{U}_{n}^{\prime}:=\left\{ U_{n}\left(I_{2n-1}\otimes(P\oplus I_{n-1})\right)U_{n}^{*}:P\text{ is a \ensuremath{n\times n} permutation matrix}\right\}.
\]

 In a \emph{ directed graph} or \emph{digraph} $\Gamma = (V,E)$, the edge set $E$ is a set of  ordered pairs $(i,j)$ with $i,j\in V$. Thus a digraph allows loops (arcs of the form $(i,i)$) and allows both $(i,j)$ and $(j,i)$ when $i\ne j$, but not multiple identical arcs between a pair of vertices.  To an arbitrary function $f: \mathbb{Z}_{n}\to{\mathbb{Z}_{n}}$ we associate a \emph{functional directed graph}, denoted by $\Gamma_{f}$, whose vertex set and arc set are respectively
\[
V\left(\Gamma_{f}\right):=\mathbb{Z}_{n},\quad E\left(\Gamma_{f}\right):=\left\{ \left(i,f\left(i\right)\right)\,:\,i\in\mathbb{Z}_{n}\right\}.
\]

Of course not every digraph arises from a function, i.e., not every digraph is functional, but the mapping $f\to \Gamma_f$ is injective.
{A digraph $\Gamma$ can  be mapped to its \emph{underlying simple graph} $G$ that has the edge $\{i,j\}$ exactly when $\Gamma$ has one or both of the arcs $(i,j)$ and $(j,i)$ and $i\ne j$.
  Similarly, $\Gamma$ can  be mapped to its \emph{underlying loop-graph} $\mathfrak G$ that has the edge $\{i,j\}$ exactly when $\Gamma$ has one or both of the arcs $(i,j)$ and $(j,i)$.
When $\Gamma_f$ is a functional digraph, we denote the underlying simple graph and loop-graph by $G_f$ and $\mathfrak G_f$, respectively.}

We define the sets of \emph{contracting functions} and \emph{non-increasing functions} as follows. 

\[\con n:=
 \left\{ h: \mathbb{Z}_{n}\to  {\mathbb{Z}_{n}}\text{ such that } h(0)=0\mbox{ and } h(i)<i\mbox{ for } i\ne 0 
\right\}.
\]

\[
\nif n:=\lsb f:\Z_n\to\Z_n \text{ such that } f\left(i\right)\le i,\,\forall\,i\in\mathbb{Z}_{n}\rsb.
\]
Note that $\con n\subset \nif n$.   Observe that while in general the mapping of a functional digraph to its underlying loop-graph is not injective, with the restriction that $f\in \nif n$ the mapping $\Gamma_f\to\mathfrak G_f$ is injective, so there is a unique association of  $\mathfrak G_f$ and $\Gamma_f$.

In Proposition \ref{p:NIFgrace} we show that $\nif n$ is in one-to-one correspondence with the set of gracefully labelled loop-graphs on $n$ vertices by associating each $f\in \nif n$ with the loop-graph having edges {$\{ f(i)+n - 1  - i, f(i)\}$ for $i = 0,\ldots, n-2$ and loop $\{f(n-1), f(n-1)\}$. Note that $\{ f(i)+n - 1  - i, f(i)\}$ is never a loop    for $i=0,\dots,n-2$. 
\begin{prop}\label{p:NIFgrace}
The set of gracefully labelled loop-graphs on $n$ vertices are in one-to-one correspondence with $\nif n$.
\end{prop}
\begin{proof}
Each $f\in\nif n$ is defined by the set of ordered pairs $\{(i,f(i)):i=0,\dots,n-1\}$. Observe that the map sending $(i,f(i))$ to the multiset $\{f(i)+n-1-i, f(i)\}$ is bijective since it is represented by
\[
\left(\begin{array}{ccc}
-1 & 1 & n-1\\
0 & 1 & 0\\
0 & 0 & 1
\end{array}\right)\left(\begin{array}{c}
i\\
f\left(i\right)\\
1
\end{array}\right)=\left(\begin{array}{c}
f\left(i\right)+n-1-i\\
f\left(i\right)\\
1
\end{array}\right).
\] 
Thus the mapping which sends $f\in \nif n$ to the loop-graph having edges $\{n - 1 + f(i) - i, f(i)\}$ for $i = 0,\ldots, n-2$ and loop $\{f(n-1), f(n-1)\}$ is bijective. 
\end{proof}

If the loop graph $\mathfrak G_f$ is gracefully labeled, it follows that the blowup construction $H_f$ is uniform and thus the Hermitian matrix $I_{2n-1}\otimes\left(A_{{\mathfrak G}_f}\oplus O_{n-1}\right)$ is $\mathfrak{U}^{\prime}_{n}$-apportionable.
We now establish a result which shows that the spectra of $\U$-apportionable matrices varies quite widely. It was recently shown in \cite{MPS21} for large $n$ that every $n$-vertex functional tree admits a $\rho$-labeling. It is known that the spectra of trees are many and varied { \cite{G84,LP73}}. Consequently, there must be a wide ranging family of spectra of matrices devised from our blowup construction of trees that are $\U$-apportionable, albeit with many repeated eigenvalues and additional zero eigenvalues. We can also invoke the stronger result from \cite{G23} that every tree admits a graceful labeling. We show that for all $f:\Z_n \rightarrow \Z_n$ subject to the iterative fixed point condition $\left|f^{\left(n-1\right)}\left(\mathbb{Z}_{n}\right)\right|=1$, the matrix $I_{2n-1}\otimes\left(\text{diag}\left(\Lambda_{\mathfrak{G}_f}\right)\oplus O_{n-1}\right)$ is $\U$-apportionable (where $\Lambda_{\mathfrak{G}_f}$ is an ordering of the spectrum of $A_{\mathfrak{G}_f}$). 

Observe that any loop-tree can be relabeled to be the loop graph of a function in $\con n$.
We have restated the next lemma from \cite{G23} in the language of loop-graphs (this is permitted because for $f\in\con n$ the mapping $\Gamma_f\to\mathfrak G_f$ is injective).

\begin{lem}\label{Lem:compo_lem}{\rm (Composition Lemma) \cite{G23}}
Let $f \in \con{n}$ be such that the length of the path whose endpoints in  $\mathfrak G_f$ are $0$ and $n-1$ is equal to the diameter of $\mathfrak G_f$ and the set $f^{-1}\left(\left\{ f\left(n-1\right)\right\} \right)$ consists of consecutive integers including $n-1$, i.e., $
f^{-1}\left(\left\{ f\left(n-1\right)\right\} \right)=\left\{ n-1,\,n-2,\,...,\,n-\left|f^{-1}\left(\left\{ f\left(n-1\right)\right\} \right)\right|\right\} 
$.
If $g \in \con{n}$ is defined from $f$ by 
\[
g\left(i\right)=\begin{cases}
\begin{array}{cc}
f^{\left(2\right)}\left(i\right) & \text{  if  }i\in f^{-1}\left(\left\{ f\left(n-1\right)\right\} \right)\\
f\left(i\right) & \text{otherwise}
\end{array},\end{cases}
\]
then
\[
\max_{\pi \in\text{S}_{n}}\left|\left\{ \left|\pi  g\pi^{-1}\left(i\right)-i\right|:i\in\mathbb{Z}_{n}\right\} \right|\le\max_{\pi \in\text{S}_{n}}\left|\left\{ \left|\pi f\pi^{ -1}\left(i\right)-i\right|:i\in\mathbb{Z}_{n}\right\} \right|.
\]
\end{lem}

\begin{thm}\label{t:Con2apport2}
If 
$f:\Z_n \rightarrow \Z_n$ subject to the iterative fixed point condition $\left|f^{\left(n-1\right)}\left(\mathbb{Z}_{n}\right)\right|=1$,
then 
\[
A=I_{2n-1}\otimes\left(A_{{\mathfrak G}(f)}\oplus O_{n-1}\right)
\]
is $\U$-apportionable.
\end{thm}

\begin{proof}  By hypothesis, there is a permutation $\pi \in S_n$ such that $\pi^{-1}f\pi\in  \con n$.  
Starting from any member of $\con{n}$, if we repeatedly perform the local iteration described in the statement of the Composition Lemma   \ref{Lem:compo_lem}, the resulting sequence of functions converges to 
 the identically zero function whose   loop-graph (a star with zero at the center) is gracefully labeled. The composition lemma asserts that the local iteration transformation never increases the maximum number of induced edge labels. Therefore the loop-graph of all members of the said sequence have graceful loop-graphs. Thus, for all functions $f\in \con{n}$ we have
\[
n=\max_{\pi\in\text{S}_{n}}\left|\left\{ \left|\pi f\pi^{ -1}\left(i\right)-i\right|:i\in\mathbb{Z}_{n}\right\} \right|.
\]
Consider the map that associates with an arbitrary $f \in \con{n}$ the $\left(2n-1\right)^{2}\times\left(2n-1\right)^{2}$ cyclic blowup matrix
\[
H_{f}=U_{n}\left(I_{2n-1}\otimes\left(A_{{\mathfrak G}(f)}\oplus O_{n-1}\right)\right)U_{n}^{*}.
\]
By Theorem \ref{t:Con2apport}  there exist a unitary matrix in the subgroup $\mathfrak{U}_{n}^{\prime}$ which apportions $H_f$.
\end{proof}

When symmetry is removed, the situation is very different. For our next result, consider the map which associates with
an arbitrary  function $f:\mathbb{Z}_{n}\to {\mathbb{Z}_{n}}$ a matrix of size $\left(2n-1\right)^{2}\times\left(2n-1\right)^{2}$
 as follows:
\[
f\mapsto T_{f}=U_{n}\left(I_{2n-1}\otimes\left(A_{f}\oplus O_{n-1}\right)\right)U_n^{*},\label{Matrix Construction 2}
\]
where $A_f \in \left\{0,1\right\} ^{n\times n}$ denotes the adjacency matrix of the functional directed graph $ \Gamma_f$ of $f$ with entries given by
\[
\lp A_{f}\rp_{i,j}=\begin{cases}
\begin{array}{cc}
1 & \text{ if }j=f\left(i\right)\\
0 & \text{otherwise}
\end{array}, & \forall\,0\le i,j<n.\end{cases}
\]
This matrix construction results in a non-Hermitian matrix unless $f$ is an involution. Observe that the subset of matrices $ \left\{ T_{f}:f:\mathbb{Z}_{n}\to {\mathbb{Z}_{n}}\right\} $ yields a matrix representation of the
transformation monoid  of functions from $\mathbb{Z}_{n}$ to ${\mathbb{Z}_{n}}$ prescribed
by the  antihomomorphism identity
\[
T_{f}\,T_{g}=T_{g\circ f}, \text{ for all }{f,g:
\Z_n\to \Z_n}
\]
 (because the adjacency matrix of a digraph acts on the right to  identify the adjacencies of a vertex).

\begin{defn}
For $A\in\Cnn$, define
\[
{\mathfrak u}(A)=\inf_{U\in {\mathfrak U}_n}\mnorm{UAU^*}.
\]
The \emph{$\mathfrak U_n$-apportionment gap} of an arbitrary matrix $A\in \mathbb{C}^{n\times n}$ is $\left|\mathfrak u(A)-\frac{\left\Vert A\right\Vert _{\text{F}}}{n}\right|$.
\end{defn}

\begin{rem}\label{r:Un-compact}
    Since the set ${\mathfrak U}_n$ is compact, the infimum is actully the minimum:
\[
{\mathfrak u}(A)=\min_{U\in {\mathfrak U}_n}\mnorm{UAU^*}.
\]
\end{rem}

\begin{thm}\label{t6.7}
If
$f\in \con n$, then
$T_{f}$ is not $\mathfrak{U}_{n}$-apportionable. 
Furthermore, the $\mathfrak{U}_{n}$-apportionement gap for every $T_f$ is
\[
\mathfrak{u}(T_{f})-\frac{\left\Vert T_{f}\right\Vert _{\text{F}}}{\left(2n-1\right)^{2}}=\frac{1}{2n-1}-\frac{\sqrt{\left(2n-1\right)n}}{\left(2n-1\right)^{2}}.
\]
\end{thm}

\begin{proof}
By Theorem \ref{t:Con2apport}, the matrix 
\[
U_{n}\left(I_{2n-1}\otimes\left(A_{f}\oplus O_{n-1}\right)\right)U_{n}^{*}+U_{n}\left(I_{2n-1}\otimes\left(\left(A_{f}-E_{0,0}\right)^{\top}\oplus O_{n-1}\right)\right)U_{n}^{*}
\]
is $\mathfrak{U}_{n}$-apportionable with apportion constant $\lp2n-1\rp^{-1}$. For any $V\in \mathfrak{U}_{n}$, the set of indices of the nonzero entries of each  $\lp 2n-1\rp \times \lp 2n-1\rp$ block indexed by $(i,j)$ of $V T_f V^*$ are identical (except possibly for some zeros created by cancellation).
Because $V$ has at its core a permutation, this set of indices is disjoint from the set of indices of the non-zero entries  in the corresponding $\lp 2n-1\rp \times \lp 2n-1\rp$  block indexed by $(i,j)$ of 
\[
VU_{n}\left(I_{2n-1}\otimes\left(\left(A_{f}-E_{0,0}\right)^{\top}\oplus O_{n-1}\right)\right)U_{n}^{*}V^*.
\]
This implies that
\begin{equation}\label{eqn: orthogonal remark}
 \left(V T_f V^{*}\right)\circ\left(VU_{n}\left(I_{2n-1}\otimes\left(\left(A_{f}-E_{0,0}\right)^{\top}\oplus O_{n-1}\right)\right)U_{n}^{*}V^{*}\right)= O_{\left(2n-1\right)^{2}}
\end{equation}
 where $\circ$ denotes the entrywise product. Therefore the action on $T_f$ by similarity transformation mediated by members of $\mathfrak{U}_{n}$ can not result in a uniform matrix. 
We conclude that $T_f$ is not $\mathfrak{U}_{n}$-apportionable as claimed.  This also implies $\mathfrak u({T_{f}})=\left(2n-1\right)^{-1}$ and
\[
\mathfrak{u}(T_{f})-\frac{\left\Vert T_{f}\right\Vert _{\text{F}}}{\left(2n-1\right)^{2}}=\frac{1}{2n-1}-\frac{\sqrt{\left(2n-1\right)n}}{\left(2n-1\right)^{2}}. \qedhere
\]
\end{proof}

As a corollary of Theorem \ref{t6.7}, for all $f,g$ lying in the semigroup $\con{n}$ the following property holds:
$\mathfrak{u}\left(T_{g}\right)=\mathfrak{u}\left(T_{f}\right).$ 

\begin{rem}
    
Observe that the three matrix summands 
\[
I_{2n-1}\otimes\left(E_{0,0}\oplus O_{n-1}\right),\, I_{2n-1}\otimes\left(\left(A_{f}-E_{0,0}\right)\oplus O_{n-1}\right) \text{ and }I_{2n-1}\otimes\left(\left(A_{f}-E_{0,0}\right)^{\top}\oplus O_{n-1}\right)
\]
are pair-wise orthogonal when viewed as members the standard (vector) inner-product space of matrices and remain so after any unitary similarity transform.

\end{rem}


\section{Application of gracefully labelled graphs to spectral inequalities}\label{s:interlacing} 

In this section we generalize the well-known  Eigenvalue Interlacing   Inequalities (involving deletion of one or more rows and columns) to a new eigenvalue interlacing inequality involving matrices obtained by zeroing the entries  associated with edges of a gracefully labeled loop-graph.  Of course there are complications for interlacing when zeroing entries is done instead of deletion even with a star (replacing deletion of row and column $j$ by zeroing row and column $j$).  Our results describe how  a combination of $n$ permutations of $n-1$ edge zeroings in a  dense arbitrarily weighted undirected graph give an interlacing bound, and the proof  uses the original Eigenvalue Interlacing Inequality.

\begin{defn}\label{d:Mk} Let $\mathfrak G$  be a loop-graph 
and recall that $A_{\mathfrak G}$ denotes its adjacency matrix. 
Recall that $C_n$ denotes the $n \times n$ cyclic shift matrix.  For $n\ge 3$ and any 
{matrix $M\in\Cnn$,  define}
\[M_{k}=M\circ\left(\left(C_n\right)^{k}\left(J_n-A_{\mathfrak G}\circ\left(J_n+I_n\right)\right)\left(C_n\right)^{-k}\right)
\] 
and denote $\spec(M_k)=\{\lambda_{k,0},\dots,\lambda_{k,n-1}\}$.
\end{defn}

We illustrate this definition in the next example.
\begin{ex}\label{e:Mk}
Let $\mathfrak G$  be the  gracefully labelled loop-graph  with edge set
\[
E\left(\mathfrak G\right)=\{ \{0,2\},\,\{0,3\},\,\{1,2\},\,\{2,2\}\} 
\]
i.e., $\mathfrak G$ is a  path on $4$ vertices with path order 3,0,2,1 and a loop at 2. Then
\[
A_{\mathfrak G}=\left(\begin{array}{rrrr}
0 & 0 & 1 & 1\\
0 & 0 & 1 & 0\\
1 & 1 & 1 & 0\\
1 & 0 & 0 & 0
\end{array}\right), \ A_{\mathfrak G}\circ\left(J_n+I_n\right)=\left(\begin{array}{rrrr}
0 & 0 & 1 & 1\\
0 & 0 & 1 & 0\\
1 & 1 & 2 & 0\\
1 & 0 & 0 & 0
\end{array}\right), \ J_n-\lp\ A_{ \mathfrak G}\circ\left(J_n+I_n\right)\rp=\left(\begin{array}{rrrr}
1 & 1 & 0 & 0\\
1 & 1 & 0 & 1\\
0 & 0 & -1 & 1\\
0 & 1 & 1 & 1
\end{array}\right)
\]
Observe that if we take
\[
M=\left(\begin{array}{rrrr}
m_{00} & m_{01} & m_{02} & m_{03} \\
m_{10} & m_{11} & m_{12} & m_{13} \\
m_{20} & m_{21} & m_{22} & m_{23} \\
m_{30} & m_{31} & m_{32} & m_{33}
\end{array}\right),
\]
then 
\[
M_0=\left(\begin{array}{rrrr}
m_{00} & m_{01} & 0 & 0 \\
m_{10} & m_{11} & 0 & m_{13} \\
0 & 0 & -m_{22} & m_{23} \\
0 & m_{31} & m_{32} & m_{33}
\end{array}\right),\qquad
M_1=\left(\begin{array}{rrrr}
m_{00} & 0 & m_{02} & m_{03} \\
0 & -m_{11} & m_{12} & 0 \\
m_{20} & m_{21} & m_{22} & 0 \\
m_{30} & 0 & 0 & m_{33}
\end{array}\right),\]
\[
M_2=\left(\begin{array}{rrrr}
-m_{00} & m_{01} & 0 & 0 \\
m_{10} & m_{11} & 0 & m_{13} \\
0 & 0 & m_{22} & m_{23} \\
0 & m_{31} & m_{32} & m_{33}
\end{array}\right),\qquad
M_3=\left(\begin{array}{rrrr}
m_{00} & 0 & m_{02} & m_{03} \\
0 & m_{11} & m_{12} & 0 \\
m_{20} & m_{21} & m_{22} & 0 \\
m_{30} & 0 & 0 & -m_{33}
\end{array}\right).
\]  
{Observe that 
\beq\label{eq:sumMj} M_0+M_1+M_2+M_{ 3}=\left(\begin{array}{rrrr}
2m_{00} & 2m_{01} & 2m_{02} & 2m_{03} \\
2m_{10} & 2m_{11} & 2m_{12} & 2m_{13} \\
2m_{20} & 2m_{21} & 2m_{22} & 2m_{23} \\
2m_{30} & 2m_{31} & 2m_{32} & 2m_{33}
\end{array}\right)=(4-2)M.\eeq
}\end{ex}

We note that with $M_k$ as defined in Definition \ref{d:Mk}, the property illustrated in \eqref{eq:sumMj} is true in general (provided $\mathfrak G$ is gracefully labelled).  The effects of the permutation similarity $\lp C_n\rp^k$ described in the proof can be observed in Example \ref{e:Mk}.

\begin{prop}\label{p:Mk}
Let $\mathfrak G$ be a gracefully labeled loop-graph of order $n\ge 4$ and let $M_k$ be as defined in Definition \ref{d:Mk}.  Then
\[
\left(n-2\right)M=\sum_{0\le k<n}M_{k}.
\]
\end{prop}
\bpf Consider the effect of the permutation similarity $\lp C_n\rp^k$: For $M=[m_{ij}]$, $\lp\lp C_n\rp^k M \lp C_n\rp^{-k}\rp_{ij}=m_{i-k,j-k}$ (with arithmetic done modulo $n$).  Because $\mathfrak G$ is gracefully labelled, the $n-1$ zeros above the diagonal in $M_0$ will land in each off-diagonal position exactly once 
 as $k$ ranges over $0,\dots,n-1$, and similarly for the zeros below the diagonal. Thus $\sum_{k=0}^{n-1}\lp M_{k}\rp_{ij}=(n-2)m_{ij}$.
For the   diagonal, observe that $\mathfrak G$ has a unique loop, so  $A_{\mathfrak G}$ has exactly one nonzero diagonal entry $m_{\ell\ell}$.  The effect of the cyclic permutation is that this one nonzero entry, which transforms  via $J_n-\lp\ A_{ \mathfrak G}\circ\left(J_n+I_n\right)\rp$ from positive to negative, hits every index once.  Thus for each $j=0,\dots,n-1$, $\sum_{k=0}^{n-1}\lp M_{k}\rp_{jj}=(n-1)m_{jj}-m_{jj}=(n-2)m_{jj}$.  \epf 

\begin{defn}\label{d:mu_k} For $M_k$ and $\lambda_{k,t}$  as defined in Definition \ref{d:Mk}, order the multiset $\{\lambda_{k,t}:t=0,\dots,n-1, k=0,\dots,n-1\}$ in nonincreasing order and denote these values by $\theta_j, j=0,\dots, n^2-1$, so that  
\[\{\theta_j, j=0,\dots, n^2-1\}=\{\lambda_{k,t}:t=0,\dots,n-1, k=0,\dots,n-1\}\mbox{ and } \theta_0\ge \theta_1\ge \dots\ge \theta_{n^2-1}.
\] 
\end{defn}

Recall that all gracefully labelled loop graphs can be constructed from nonincreasing functions (Proposition \ref{p:NIFgrace}). 
\begin{thm}\label{t:spec-ineq} Let $\mathfrak G$ be a gracefully labelled loop-graph of order $n\ge 4$, let $M$ be an $n\x n$ Hermitian matrix, and let $\theta_j$ be as defined in Definition \ref{d:mu_k}.   For $\ell=0,\dots,n-1$,
\[\frac n{n-2}\theta_{\ell+(n^{2}-n)}\le \lam_{\ell}(M)\le \frac n{n-2}\theta_{\ell}.\]
\end{thm}
\bpf 
Since $M$ is Hermitian, each $M_k$ is a Hermitian matrix. Thus by the Spectral Theorem, each matrix $M_k$  admits a spectral decomposition of the form
\beq\label{eq:decomp1}
M_{k}=U_{k}\diag\left(\Lambda_{k}\right)U_{k}^{*},
\eeq
where $U_{k}$ is a real unitary matrix,   $\Lambda_{k}=\lp\lambda_{k,0},\dots,\lambda_{k,n-1}\rp $ and $\spec(M_k)=\{\lambda_{k,0},\dots,\lambda_{k,n-1}\}\subset \R$.
Thus 
\beq\label{eq:decomp2} \lp M_k\rp_{ij}=\sum_{0\le t<n}
\lp\lp U_k\rp_{i,t}\sqrt{\lambda_{k,t}} \rp
\lp\sqrt{\lambda_{k,t}}\,  \lp \overline{ U_k}\rp_{j,t} \rp.\eeq
From Proposition \ref{p:Mk}, 
\beq\label{eq:decomp3}
m_{ij}=\sum_{0\le k<n}\sum_{0\le t<n}\left(\frac{\lp U_k\rp_{i,t}}{\sqrt{n}}\,\sqrt{\frac{n}{n-2}\,\lambda_{k,t}}\right)\, \left(\sqrt{\frac{n}{n-2}\,\lambda_{k,t}}\,\frac{\lp \overline{ U_k}\rp_{j,t}}{\sqrt{n}}\right),\ \forall\:0\le i,j<n.
\eeq
Reversing the process used to go from \eqref{eq:decomp1} to \eqref{eq:decomp2}, we view the entry-wise equality in \eqref{eq:decomp3} as expressing the product of  three matrices. The first matrix is $n\times n^2$ matrix $\hat U$ defined as follows: The $i$-th row of $\hat U$ is obtained by concatenating row $i$ of the $n$ matrices $U_k$ for $k=0,\dots,n-1$. The second matrix is the $n^2\x n^2$ diagonal matrix
\[
\Lambda=\bigoplus_{k=0}^{n-1}\text{diag}\left(\frac n{n-2}{\Lambda}_{k}\right).
\]
  Finally the third matrix is the Hermitian adjoint of the first matrix $\hat U$.  Observe that the rows of $\hat U$ are orthonormal. By extending the rows of $\hat U$ to an orthonormal basis for $\R^{n^2}$ and applying  Gram-Schmidt, we can expand   $\hat U$  to a unitary matrix $U$ of size $n^2 \times n^2$.  
  
 Then $M$ is the $\{0,\dots,n-1\}$ principal submatrix of the $n^2 \x n^2$ matrix $U\Lambda U^*$, i.e.,
\[
U\Lambda U^{*}=\mtx{M & B_{0,1}\\{B_{0,1}}^* & B_{1,1}}
\]
for some matrices $B_{0,1}$ and $B_{1,1}$.
Then by the Eigenvalue Interlacing Theorem \cite[Theorem 8.10]{Zhang},
\[\lambda_{\ell+(n^{2}-n)}\left(\left[\begin{array}{cc}
M & B_{0,1}\\
B_{0,1}^{*} & B_{1,1}
\end{array}\right]\right)\le\lambda_{\ell}\left(M\right)\le\lambda_{\ell}\left(\left[\begin{array}{cc}
M & B_{0,1}\\
B_{0,1}^{*} & B_{1,1}
\end{array}\right]\right)\] 
for $\ell=1,\dots, n$.  
Observe that  \[\frac n{n-2} \theta_{\ell}=\lambda_{\ell}\lp \Lambda \rp =\lambda_{\ell}\left(\left[\begin{array}{cc}
M & B_{0,1}\\
B_{0,1}^{*} & B_{1,1}
\end{array}\right]\right). \qedhere\]
\epf

\section{Spectra and Jordan canonical forms 
of GL-apportionable matrices}\label{s:spectra}

 We now shift focus from unitary apportionability to general apportionability. In this section we study the question of
what multisets of $n$ complex numbers can and cannot be realized as spectra  or Jordan canonical forms  of uniform $n\x n$  matrices.
Being  the spectrum or Jordan canonical form of an apportionable matrix is equivalent to being the spectrum or  Jordan canonical form of a uniform matrix since similarity by matrices in $\gl n$ is allowed for GL-apportionability. 
We begin with some elementary observations and then focus on the case of $2\x 2$ matrices.

For every $\lambda\in \C$, there is a uniform matrix  $B\in\Cnn$ with $\rank B=1$ and $\spec(B)=\{\lambda,0,\dots,0\}$ by Theorem \ref{t:u-apport-exist}. 
We can scale the spectrum of a uniform  matrix: If $\Lambda=\spec(B)$ for a uniform  matrix $B$, then  for $b\in\C$,  $b\Lambda=\spec(bB)$ and $bB$ is uniform.

Kronecker products can be used to construct bigger uniform matrices from smaller uniform matrices, and thus expand the set of spectra that we are able to realize with uniform matrices.  For two multisets $S=\{s_0,\dots,s_{n-1}\}$ and $T=\{t_0,\dots,t_{m-1}\}$, define $ST=\{s_kt_j:k=0,\dots,n-1, j=0,\dots,{m-1}\}$.

\begin{rem}
Suppose $B_n$ and $B_m$ are uniform $n\x n$ and $m\x m$ matrices, with spectra $\Lambda_n$ and $\Lambda_m$. Then $B_n\otimes B_m$ is uniform.  From known properties of Kronecker products \cite[Theorem 4.8]{Zhang},  $\spec(B_n \otimes B_m)=\Lambda_n \Lambda_m $. By using $B_n=F_n$ (recall that $F_n$ denotes the $n\x n$ DFT matrix) and $B_m$ as any uniform matrix with spectrum $\Lambda$, we conclude that if $\Lambda$ is the spectrum of a uniform matrix, then so is $\cup_{j=1}^n \omega^j \Lambda$ where $\omega $ is a primitive $n$th root of unity.
\end{rem}

\begin{prop}
Let $B$ be a uniform $n\x n$ matrix. 
Then there is a uniform matrix with spectrum $\spec(B) \cup\{0^{((r-1)n)}\}$ for any positive integer $r$.
\end{prop}
\bpf Let $E_{0,0}$ be the $r\x r$  matrix with $(0,0)$-entry equal to 1 and all other entries 0.  Consider $B'=E_{00}\otimes B=B\oplus O_{(r-1)n}$. Then $\spec(B')= \spec(B) \cup\{0^{((r-1)n)}\}$.   Define $U=F_r\otimes I_n$ and observe that $U$ is unitary.  Furthermore,
\[UB'U^*=(F_r\otimes I_n)(E_{00}\otimes B)({F_r^*}\otimes I_n)=(F_r E_{00} F_r^*)\otimes B=J_r\otimes B\]
is uniform. \epf

Since extending the spectrum with blocks of zeros preserves apportionability, it is natural to ask whether we can add just one zero and preserve apportionability.
\begin{quest}
If  $B$ is uniform, is $\spec(B) \cup \{0\}$ the spectrum of a uniform matrix?
\end{quest}

 Next we study the question of
what multisets of $2$ complex numbers cannot be realized as spectra of uniform $2\x 2$  matrices, and thus of apportionable matrices.  
Very few real spectra are realizable by uniform matrices.

\begin{thm}\label{constant-spec}
For any $n\ge 2$, the spectrum $\{0^{(n)}\}$ can be realized by a nonzero  uniform matrix.   For any nonzero $\lam \in\C$, the spectrum $\{\lam^{(2)}\}$ cannot be realized by a uniform matrix.
\end{thm}
\bpf Since $\rank E_{0,1}=1$ where $E_{0,1}$ is $n\x n$, Theorem \ref{t:u-apport-exist} shows $\{0^{(n)}\}$ can be realized as the spectrum of a nonzero uniform matrix.

Suppose that  $B\in\C^{2\x 2}$ is a uniform matrix with spectrum $\{\lam^{(2)}\}$.  We show that $\lam=0$. There exists a matrix $M$ such that  $M^{-1}BM$ is in Jordan canonical form. Since $B$ is uniform and $M\lam I_2M^{-1}=\lam I_2$, if the Jordan canonical form of $B$ is $\lam I_2$, then $\lam=0$.   So assume the Jordan canonical form of $B$ is  $\lam I_2+N$ where $N=\mtx{0 & 1\\0 & 0}$.  Then $M\lp \lam I_2+N\rp M^{-1}=\lam I_2+MNM^{-1}$.  Without loss of generality, $\det M$=1 and $M=\mtx{x & y\\z &\frac {1+yz}x}$.    Then 
\[B=M\lp \lam I_2+N\rp M^{-1}=\mtx{\lam & 0\\0 & \lam}+\mtx{-xz  & x^2\\-z^2 & xz}=\mtx{-xz  +\lam& x^2\\-z^2 & xz+\lam}\]
is uniform, so $|x|=|z|$.  If $x=0$, then $\lam=0$, so assume $x\ne 0$.    Now uniformity implies 
$|xz+\lam|=|-xz+\lam|=|x^2|=|xz|.$
Let $\delta=xz$.  Then $|xz+\lam|^2=\delta\ol{\delta}+\delta\ol{\lam}+\lam\ol{\delta}+\lam\ol{\lam}$, $|-xz+\lam|^2=\delta\ol{\delta}-\delta\ol{\lam}-\lam\ol{\delta}+\lam\ol{\lam}$, and $|xz|^2=\delta\ol{\delta}$.  Thus $\delta\ol{\lam}+\lam\ol{\delta}=0$. So $\delta\ol{\delta}=\delta\ol{\delta}+\lam \ol{\lam}$, which implies $\lam=0$.
\epf

The next result is established by computation.
\begin{lem}\label{lem:2x2MAMinv}
Let $c\in \C$, $D=\diag(1,c)$ and  $M=\mtx{x & y\\z &\frac {1+yz}x}$ Then \[B=MDM^{-1}=\mtx{1+(1-c)yz & -(1-c)xy\\(1-c)(1+yz)\frac z x & c-(1-c)yz}. \]
\end{lem}

\begin{thm}\label{2x2-real-spec}
For a real number $r$, the spectrum $\{1,r\}$ can be realized by a uniform matrix if and only if $r=0$ or $r=-1$.
\end{thm}
\bpf The spectra $\{1,0\}$ and $\{1,-1\}$ are realized by $\frac 1 2 J_2$ and $\frac 1 {\sqrt 2}H_2$, respectively (recall $H_2$ is a $2\x 2$ Hadamard matrix).

The eigenvalues of a uniform $2\x 2$ matrix are distinct unless both are zero 
 by Theorem \ref{constant-spec}.   Let $r\in \R$ and $D=\diag(1,r)$.  Assume $D$ is apportionable.  Then we may assume the apportioning matrix $M$ has the form in Lemma \ref{lem:2x2MAMinv}, so
$B=\mtx{1+(1-r)yz & -(1-r)xy\\(1-r)(1+yz)\frac z x & r-(1-r)yz}$ is uniform.
Let $yz=a+b \ii$ with $a,b\in\R$.  Compare the absolute values of the  $(0,0)$ and $(1,1)$ entries of $B$:
\bea |1+(1-r)(a+b\ii)|&=&|r-(1-r)(a+b \ii)|.\\
|(1+(1-r)a) +((1-r)b) \ii|&=&|(r-(1-r)a) +((1-r)b) \ii|.\\
|1+(1-r)a|&=&|r-(1-r)a|.
\eea

Thus $1+(1-r)a=r-(1-r)a$ or $-(1+(1-r)a)=r-(1-r)a$. Since $-(1+(1-r)a)=r-(1-r)a$ implies $r=-1$ and we have seen that $r=-1$ can be realized, assume $1+(1-r)a=r-(1-r)a$.  Thus  $a=-\frac 1 2$.  
Since $B$ is uniform, the absolute value of product of the off-diagonal entries equals the absolute value of product of the diagonal entries.  Since $1+yz=-\ol{yz}$, the product of the off diagonal entries is $-(1-r)^2yz(1+yz)=(1-r)^2|yz|^2=(1-r)^2(\frac 1 4+b^2)$.
The square of the absolute value of each entry must be $|1+(1-r)yz|^2=|1+(1-r)(-\frac 1 2+b \ii)|^2=|\frac{1+r}2+(1-r)b\ii|=\frac{(1+r)^2}4+ (1-r)^2b^2$.  Thus
\bea (1-r)^2(\frac 1 4+b^2)&=& \frac{(1+r)^2}4+ (1-r)^2b^2\\
\frac{(1-r)^2} 4&=& \frac{(1+r)^2}4\\
0&=&r.
\eea
Thus the uniformity of  $B$  implies $r=-1$ or $r=0$.
\epf

  It is immediate from the previous theorem that  two nonzero eigenvalues of a $2\x 2$ uniform matrix may or may not have the same magnitude.  This is also  illustrated in the next two examples.

\begin{ex}\label{ex1I}  For $B=\frac 1 2\mtx{1+\ii & -1+\ii\\-1+\ii& 1+\ii}$,  $\spec(B)=\{1,\ii\}$.
\end{ex}

\begin{ex}\label{diff-mag-evals}  For $A=\mtx{1 & 1\\1& \frac{1+\sqrt 3 \ii}2}$,  the (approximate) decimal values of the eigenvalues of $A$ are $1.69244 + 0.318148 \ii$ and $ -0.19244 + 0.547877 \ii$.
\end{ex}


\section{Finding an apportioning matrix $M$ and constant  $\kappa$}\label{s:M-kappa}
In this section we discuss how to find GL-apportioning matrices. We begin with a simple $2\times 2$ example that illustrates a matrix can have infinitely many GL-apportionment constants each of which can be obtained from infinitely many apportioning matrices.

\begin{ex}\label{e:app-const-not-unique}
Let 
$A = \left[\begin{array}{cc}
2 & 0\\
0 & 0
\end{array}\right]$ 
and let 
$M =
\left[\begin{array}{cc}
w & x\\
y & z
\end{array}\right]$ 
be nonsingular (so $wz\ne xy$). Then
\[
MAM^{-1} = \frac{2}{wz - xy} 
\left[\begin{array}{cc}
wz & -wx\\
yz & -xy
\end{array}\right].
\]
Thus the  matrix $MAM^{-1}$ is uniform if and only if $|x|=|z|$ and $|w|=|y|$. Observe that $|wz - xy| \leq 2|wz|$ and so each apportionment constant for $A$ must be at least $1$. Let $a,\theta\in \R$ such that $a\not=0$ and $0 <\theta < 2\pi$, and let
\[
M = 
\left[\begin{array}{cc}
a & a^{-1} \\ 
a & a^{-1}e^{i\theta}
\end{array}\right].
\]
Then $M$ is nonsingular and apportions $A$ with apportionment constant $\kappa = \big|\sin\big(\frac\theta 2\big)\big|^{-1}$. Notice that $\kappa = 1$ for $\theta = \pi$, and that $\kappa$ can be made arbitrarily large for a sufficiently small choice of $\theta$. Thus $[1,\infty)$ is the set of apportionment constants for $A$.
\end{ex}

Example \ref{e:app-const-not-unique} utilizes ad hoc methods to solve for apportioning matrices of a small and curated matrix. It may seem rather hopeless to find apportioning matrices in a more general setting. The search for apportioning matrices can be simplified with the following  proposition. Let $\mvec(A)$ denote the vectorization of the matrix $A$.
Recall that $\circ$ denotes the entrywise product.

\begin{prop}\label{prop:additive_apport_mat}
Let $A\in\C^{n\times n}$ and $M\in \gl n$. Let $\bv = \mvec\lp (MAM^{-1})\circ (\overline{MAM^{-1}})  \rp$ and let $F$ be the $n^2\times n^2$ DFT matrix. Then $M$ apportions $A$ if and only if $F \bv \in \spn(\be_0)$.
\end{prop}
\begin{proof}
Suppose that $M$ apportions $A$ and let $\kappa$ be the apportioning constant for $M$. Then $(MAM^{-1})\circ (\overline{MAM^{-1}}) = \kappa^2 J$ and so $F\bv = n \kappa^2 \be_0$.

Now suppose that $F \bv \in \spn(\be_0)$. Then $\bv = c \bone$ for some $c\in\R$ (by construction $\bv\in \R^{n^2}$) and hence $(MAM^{-1})\circ (\overline{MAM^{-1}}) = c J$. Thus $MAM^{-1}$ is uniform and so $M$ apportions $A$.
\end{proof}

Proposition \ref{prop:additive_apport_mat} can be used to solve for apportioning matrices by generating a system of $n^2-1$ equations in the entries of $M$. Note that $F$ can be replaced with any unitary matrix whose first row is a multiple of $\bone$. 

We revisit Example \ref{e:app-const-not-unique} to illustrate how to apply Proposition \ref{prop:additive_apport_mat}.

\begin{ex}
Let $A$ and $M$ be the same as in Example \ref{e:app-const-not-unique}. We may assume, without loss of generality, that $\det(M) = wz - xy = 2$. Then
\[
(MAM^{-1})\circ (\overline{MAM^{-1}}) = 
\left[\begin{array}{cc}
|wz|^2 & |wx|^2\\
|yz|^2 & |xy|^2
\end{array}\right].
\]
By Proposition \ref{prop:additive_apport_mat}
\[
\begin{aligned}
|wz|^2 + \ii|yz|^2 - |wx|^2 - \ii|xy|^2 &= 0,\\
|wz|^2 - |yz|^2 + |wx|^2 - |xy|^2 &= 0,\\
|wz|^2 - \ii|yz|^2 - |wx|^2 + \ii|xy|^2 &= 0.
\end{aligned}
\]
This system of equations can be reduced to $|x| = |z|$ and $|w| = |y|$.\end{ex}

\begin{rem}
Suppose that $A\in \mathbb{C}^{n\times n}$ is $\U$-apportionable. Then the entries of a unitary matrix $U$ that apportions $A$ can be determined by Proposition \ref{prop:additive_apport_mat} along with the system of equations resulting from $UU^* = I$.
\end{rem}

Note that $A^{\circ^{-1}}$ means the entrywise inverse of a  matrix $A$ because $\circ$ is the entrywise product.

\begin{thm}\label{apport_decomp}
Let $A\in\mathbb{C}^{n\times n}$ be nonzero and apportionable with an apportionment constant $\kappa\ge 0$. Then there
exists an $M\in \gl n$ such that
\[\label{apportion_decomposition}
A=\kappa^2\,M^{-1}\left(\overline{MAM^{-1}}\right)^{\circ^{-1}}M.
\]
\end{thm}

\begin{proof}
Since $A$ is apportionable with apportionment constant $\kappa$, there exists an $M\in \gl n$ such that $B = MAM^{-1}$ is uniform and $\kappa = \|B\|_{\max}$. Observe that $
B\circ \overline{B} =\kappa^2 J.
$
Since $\overline{B}$ has no zero entries, $A = \kappa^2\,M^{-1}\,\overline{B}^{\circ^{-1}}M$, as claimed.
\end{proof}

\begin{quest}\label{q:realize2} When $A\in \mathbb{C}^{n\times n}$ is not apportionable how do we find and certify the matrix $M$
that achieves  the infimum, 
$
\inf_{M\in\gl n}\left\Vert MAM^{-1}\right\Vert _{\text{max}}
$?
\end{quest}

\section{Concluding remarks}\label{s:conclude} 
We have included open questions throughout when relevant to the material discussed.  In this section we list some additional open questions.  

We begin  with questions related to how `common' apportionable matrices are.
For context, recall that set of matrices that cannot be diagonalized is of measure zero (because an eigenvalue must be repeated).  What about apportionability?  It was shown in Proposition \ref{c:notU} that the set of matrices that are not $\U$-apportionable has positive measure.

\begin{quest}  Is the set of $\U$-apportionable matrices of measure zero or positive measure? 
\end{quest}

\begin{quest} Is the set of  matrices that are not GL-apportionable of measure zero or positive measure?  Is the set of GL-apportionable matrices of measure zero or positive measure? 
\end{quest}

There are numerous ways to measure closeness to apportionability.  Section \ref{s:u-close} contains results about one such measure for $\U$-appotionability,  $u(A)=\min_{U\in \U(n)}\mnorm{UAU^*}$.  Here we mention other possibilities.
\begin{defn}\label{d:rat}
For $A\in\Cnn$ with no zero entries, define the \emph{uniformity ratio} to be $ur(A)=\frac {\max\{|a_{ij}|\}}{\min\{|a_{ij}|\}}$; if there are both zero and nonzero entries in $A$, then $ur(A)=\infty$.
Define the \emph{unitary apportionability ratio} of $A\ne O$ to be
\[
uar(A)=\inf_{U\in \U(n)}ur(UAU^*).\]
\end{defn}

Let $A\in\Cnn$ and $A\ne O$.  Observe that $A$ is $\U$-apportionable if and only if $uar(A)=1$.
A unitary matrix $U$ obtained from a random $n\x n$ matrix via orthonormalization of the columns will have the property that $UAU^*$ has no zero entries and thus $ur(A)<\infty$.

\appendix

\section{Recovery Lemma}\label{appdx}

 The Composition Lemma, which is proved in \cite{G23}, is applied in Section \ref{s:graph-decomp}.  It relies on the Recovery Lemma.  Here we provide a proof of the Recovery Lemma.

\begin{rem}\label{r:Gf}
Let $f: \mathbb{Z}_{n}\to{\mathbb{Z}_{n}}$ be a function. Recall that  the functional digraph $\Gamma_{f}$ associated with $F$ has $
V\left(\Gamma_{f}\right)=\mathbb{Z}_{n}$ and $E\left(\Gamma_{f}\right)=\left\{ \left(i,f\left(i\right)\right)\,:\,i\in\mathbb{Z}_{n}\right\}$. Each vertex in  $\Gamma_f$ has out degree one.  A fixed point of $f$ corresponds to a loop in $\Gamma_f$.  Note that $f$ can be determined from $\Gamma_f$ (but not always from  the underlying simple graph $G_f$).  If  $G_f$ is connected, then $f$ has at most one fixed point, because $n-1$ non-loop arcs are needed.  If  $G_f$ is connected and $f$ has a fixed point, then $\Gamma_f$ does not have any cycles except the loop at the fixed point.
If $f$ has a fixed point and $G_f$ is connected, then {the fixed point} and the edges of $G_f$ uniquely determine  $f$ and $\Gamma_f$ :   Let $u$ be the unique fixed point and let initially define $X=\{u\}$;  $X$ is the set of vertices $x$ for which $f(x)$ is determined.  If $v$ is a  neighbor of $x\in X$, then $f(v)=x$ and the arc is $(v,x)$, since each vertex of $\Gamma_f$ has out degree one. So now $X:=N[X]$.  Repeat the this neighborhood step until $X=\Z_n$.
\end{rem}  

\begin{defn}
For a function $ f:\mathbb{Z}_{n}\to {\mathbb{Z}_{n}}$, define the \emph{edge-labeling}  polynomial of $f$ to be \[
p_f\lp x_0,\dots,x_{n-1}\rp=\prod_{0\le i<j<n}\left(\left(x_{f\left(j\right)}-x_{j}\right)^{2}-\left(x_{f\left(i\right)}-x_{i}\right)^{2}\right),
\]
\end{defn}

We illustrate why this is called the edge-labeling polynomial in the next example.
\begin{ex}\label{ex:P4}
Let $f:\Z_4\to \Z_4$ and suppose that $G_f = P_4$ as depicted in  Figure \ref{fig1}.

\begin{figure}[!h]
\begin{center}
\begin{tikzpicture}
\node (3) at (-2,0) {};
\node (2) at (0,0) {};
\node (1) at (2,0) {};
\node (0) at (4,0) {};
\draw[fill=black] (-2,0) circle (3pt);
\draw[fill=black] (0,0) circle (3pt);
\draw[fill=black] (2,0) circle (3pt);
\draw[fill=black] (4,0) circle (3pt);
\node at (-2,0.5) {$3$};
\node at (0,0.5) {$2$};
\node at (2,0.5) {$1$};
\node at (4,0.5) {$0$};
\node at (-3,0) {$G_f =$};
\draw (3) edge[thick,-] (2);
\draw (2) edge[thick,-] (1) edge[thick,-] (0);
\end{tikzpicture}
\end{center}
    \caption{The graph $G_f$. \label{fig1}}
\end{figure}
Assume $0$ is the unique fixed point of $f$. Then $f(0) = 0$ and $f(i) = i-1$ for $0< i \leq 3$. The functional directed graph $\Gamma_f$ is shown in Figure \ref{fig2} below. 
   \begin{figure}[!h]
\begin{center}
\begin{tikzpicture}
\node (3) at (-2,0) {};
\node (2) at (0,0) {};
\node (1) at (2,0) {};
\node (0) at (4,0) {};
\draw[fill=black] (-2,0) circle (3pt);
\draw[fill=black] (0,0) circle (3pt);
\draw[fill=black] (2,0) circle (3pt);
\draw[fill=black] (4,0) circle (3pt);
\node at (-2,0.5) {$3$};
\node at (0,0.5) {$2$};
\node at (2,0.5) {$1$};
\node at (4,0.5) {$0$};
\node at (-3,0) {$\Gamma_f =$};
\draw (3) edge[thick,->] (2);
\draw (2) edge[thick,->] (1) edge[thick,->] (0);
\draw (0) edge[thick,->,out=45,in=315,looseness=10] (0);

\end{tikzpicture}
\end{center} 
    \caption{The graph $\Gamma_f$. \label{fig2}}
\end{figure}

In order to determine $p_f$ note that $f(0)=0$, and so $(x_{f(j)}-x_j)^2-(x_{f(0)}-x_0)^2 =(x_{f(j)}-x_j)^2$ for $j>0$.  Thus 
\begin{linenomath}
{\small
\begin{align*}
p_f(x_0,x_1,x_2,x_3) &= \prod_{0\leq i < j < 4}\big((x_{f(j)} - x_j)^2 - (x_{f(i)} - x_i)^2 \big)\\
&= { \lp\prod_{j=1}^4(x_{j-1} - x_j)^2 \rp 
\lp (x_1 - x_2)^2 -(x_0 - x_1)^2\rp 
\lp (x_2 - x_3)^2 -(x_0 - x_1)^2\rp 
\lp(x_2 - x_3)^2 - (x_1 - x_2)^2\rp} 
\end{align*}}
\end{linenomath}
\end{ex}

\begin{obs}
The edge-labeling polynomial $p_f\lp x_0,\dots,x_{n-1}\rp$ is not identically zero if and only if $f$ has at most one fixed point 
and $\Gamma_f$ has no $2$-cycles. 
\end{obs}

The next result gives an algorithm for recovering $G_f$ from $p_f$ when $p_f$ is not identically zero and $f$ has a fixed point.  

\begin{lem}[Recovery Lemma]\label{t:fpf-inject}
 Suppose the edge-labeling polynomial $p_f\lp x_0,\dots,x_{n-1}\rp$ is defined from some function $f:\Z_n\to \Z_n$ and $p_f$ is not identically zero.  It can be determined from $p_f$ whether or not $f$ has a (necessarily unique) fixed point. If $f$ has a  fixed point,  then $G_f$ can be determined from $p_f$.  If  $f$ has a fixed point and $G_f$   is connected, then $f$ and $\Gamma_f$ can be determined from $p_f$ and the fixed point.  Let $S$ denote the set of functions $ f:\mathbb{Z}_{n}\to {\mathbb{Z}_{n}}$ such that $f$ has a unique fixed point  $0$  and $G_f$ is connected.   The function from $S$ to $\mathbb{Q}\left[x_{0},\cdots,x_{n-1}\right]$ that assigns $p_f$ to
$f$ is injective.
\end{lem}
\bpf
We show  that each factor in a factorization of $p_f$  is   a quadrinomial (a linear combination of  exactly four distinct variables),
a trinomial (a linear combination of  exactly three distinct variables), or a binomial (a linear combination of of exactly two distinct variables), and analyze how each can occur.  

We factor 
\[\lp\left(x_{f\left(j\right)}-x_{j}\right)^{2}-\left(x_{f\left(i\right)}-x_{i}\right)^{2}\rp=\left(x_{f\left(j\right)}-x_{j}+x_{f\left(i\right)}-x_{i}\right)\left(x_{f\left(j\right)}-x_{j}-x_{f\left(i\right)}+x_{i}\right)\]
A factor $x_{f\left(j\right)}-x_{j}-x_{f\left(i\right)}+x_{i}$ or $x_{f\left(j\right)}-x_{j}+x_{f\left(i\right)}-x_{i}$ has the form $a+b-c-d$; it is is a quadrinomial  with $a,b,c,d$ distinct if and only if $|\{a,b,c,d\}|=4$, i.e., 
$
\left|\left\{ x_{f\left(j\right)},x_{j},x_{f\left(i\right)},x_{i}\right\} \right|=4
$. In this case both $x_{f\left(j\right)}-x_{j}-x_{f\left(i\right)}+x_{i}$ and $x_{f\left(j\right)}-x_{j}+x_{f\left(i\right)}-x_{i}$ are quadrinomials.

 The expression $a+b-c-d$ collapses to a binomial if $|\{a,b\}\cap \{c,d\}|=1$ (note that $|\{a,b\}\cap \{c,d\}|=2$ is impossible since $p_f$ is not identically zero).  Notice that $a+b-c-d$ occurs in two forms in $p_f$: 
$\{a,b\}=\{ x_{f(j)}, x_{f(i)}\}, \{c,d\}=\{x_j,x_i\}$ or $\{a,b\}=\{ x_{f(j)}, x_i\}, \{c,d\}=\{x_j,x_{f(i)}\}$.  First consider the case that $f$ has a (unique) fixed point $u$.  Then for each $j\ne u$ we obtain two copies of the binomial $f\lp x_j\rp -x_j$ from $\pm \lp\lp f\lp x_j\rp -x_j\rp^2-\lp f\lp x_u\rp -x_u\rp^2\rp=\pm \lp f\lp x_j\rp -x_j\rp^2$ with $+$ if $j>u$ and $-$ otherwise.

Now assume neither $i$ nor $j$ is a fixed point.  A binomial-trinomial pair of factors arises from\break$\left(x_{f\left(j\right)}-x_{j}+x_{f\left(i\right)}-x_{i}\right)\left(x_{f\left(j\right)}-x_{j}-x_{f\left(i\right)}+x_{i}\right)$ when $\{a,b\}=\{ x_{f(j)}, x_{f(i)}\}, \{c,d\}=\{x_j,x_i\}$, and $j=f(i)$ or $i=f(j)$.  Without loss of generality, we choose $j=f(i)$. 
This produces
\[
\pm  \lp  x_{f\left(f(i)\right)}-x_i\right)\left(x_{f\left(f(i)\right)}+x_{i}-2x_{f\left(i\right)}\rp.
\]  
 Similarly, a binomial-trinomial pair of factors arises when $\{a,b\}=\{ x_{f(j)}, x_i\}, \{c,d\}=\{x_j,x_{f(i)}\}$, which implies $f\left(i\right)=f\left(j\right)$. Setting $i<j$, this  produces
\[
\left(2x_{f\left(j\right)}-x_{j}-x_{i}\rp \lp x_i-x_{j}\right).
\]

We have now described all possible ways binomial factors can occur in $p_f$.  Furthermore,
a trinomial factor of $p_f$ can only occur in a binomial-trinomial pair.  Observe that in each binomial-trinomial pair, the trinomial has the form $\pm (2r-s-t)$ and the associated binomial is of the form $(s-t)$. 

We now take a given polynomial $p_f$  that is not identically zero, with no information about $f$ except that $f:\Z_n\to \Z_n$ is a function.  Define 
$h\lp x_0,\dots,x_{n-1}\rp$ to be the product of all the binomials that occur in binomial-trinomial pairs.  That is, $s-t$ is a factor of $h$ if and only if  $2r-s-t$ is a factor of $p_f$ for some $r$.  Now define \[q\lp x_0,\dots,x_{n-1}\rp=\frac {p_f\lp x_0,\dots,x_{n-1}\rp}{h\lp x_0,\dots,x_{n-1}\rp},\] which is a polynomial. Then $q$ has no binomial factors if and only if $f$ does not have a fixed point.  Otherwise, $q$ has $2(n-1)$ binomial factors, which occur in pairs: $(x_k-x_\ell)^2$.  Then  $E(G_f)=\{k \ell:(x_k-x_\ell)^2 \mbox{ is a factor of }q\}$.  The remaining two statements now follow from knowing $G_f$ by Remark \ref{r:Gf}. 
\epf


\section*{Acknowledgements}
We than Steve Kirkland and  Hermie Monterde for bringing the connections between apportionment and instantaneous uniform mixing of quantum walks to our attention. We thank the reviewers for helpful comments   that have improved the exposition of paper. \\
 The research of B. Curtis was partially supported by NSF grant 1839918.\\
 The research of Edinah K. Gnang was partially supported by Technical Information Center (DTIC) under award number FA8075-18-D-0001/0015.

\end{document}